\documentstyle[amscd, oneside]{amsart}
\font\tencmsy=cmsy10
\font\sevencmsy=cmsy7
\font\fivecmsy=cmsy5
\newfam\cmsyfam
\textfont\cmsyfam=\tencmsy
\scriptfont\cmsyfam=\sevencmsy
\scriptscriptfont\cmsyfam=\fivecmsy

\font\tenmsb=msbm10
\font\sevenmsb=msbm7
\font\fivemsb=msbm5
\newfam\msbfam
\textfont\msbfam=\tenmsb
\scriptfont\msbfam=\sevenmsb
\scriptscriptfont\msbfam=\fivemsb

\let\Bbb\Bbbb

\font\teneufm=eufm10
\font\seveneufm=eufm7
\font\fiveeufm=eufm5
\newfam\eufmfam
\textfont\eufmfam=\teneufm
\scriptfont\eufmfam=\seveneufm
\scriptscriptfont\eufmfam=\fiveeufm

\let\goth\frak

\def\gg{\goth g}

\def\gh{\goth h}
\def\gk{\goth k}

\def\gb{\goth b}
\def\bg{\goth b}
\def\gp{\goth p}

\def\Z{\Bbb Z}
\def\C{\Bbb C}
\def\R{\Bbb R}
\def\Q{\Bbb Q}
\def\sv{{\sl v}}

\def\dlim{\varinjlim}
\def\lim{{\rm lim}}

\def\sgn{{\rm sgn}}
\def\id{{\rm id}}

\def\supp{{\rm supp}}

\def\dim{{\rm dim}}

\def\rk{{\rm rk\,}}
\def\int{{\rm int}}
\def\ss{{\rm ss}}
\def\ch{{\rm ch}}

\def\cplus{\hbox{$\subset${\raise1.05pt\hbox{\kern-0.65em 
${\scriptscriptstyle +}$}}\ }}
\def\bcplus{\hbox{$\supset${\raise1.05pt\hbox{\kern-.65em 
${\scriptscriptstyle +}$}}\ }}

\def\nek{\text{\hbox{$\simeq$ \kern-.95em \hbox{$/$ \kern.05em}}}}
\def\vep{\varepsilon}
\def\opp{\operatornamewithlimits{\oplus}}

\newtheorem{Pro}{Proposition}
\newtheorem{The}{Theorem}
\newtheorem{Lem}{Lemma}

\title  {Weight modules of direct limit Lie algebras}
\author {Ivan Dimitrov and Ivan Penkov}
\address {I.D.: \quad
Department of Mathematics, University of California at
Riverside, Riverside, CA 92521, USA}
\email{dimitrov@@math.ucr.edu}
\address {I.P.: \quad
Department of Mathematics, University of California at
Riverside, Riverside, CA 92521, USA}
\email{penkov@@math.ucr.edu}
\date{\today}

\begin{document}

\begin{abstract}
In this article we initiate a systematic study of irreducible weight modules 
over direct limits of reductive Lie algebras, and in particular over 
the simple Lie algebras $A(\infty)$, $B(\infty)$, 
$C(\infty)$ and $D(\infty)$. Our main tool is the shadow method introduced
recently in \cite{DMP}.
The integrable irreducible modules are an important particular class
and  we give an explicit parametrization of
the finite integrable  modules which are  
analogues of finite-dimensional irreducible modules over reductive 
Lie algebras.
We then introduce the more general class of pseudo highest weight modules. 
Our most general result is the description of the support of any irreducible weight 
module.

Key words (1991 MSC): Primary 17B10.
\end{abstract}
\maketitle

\section*{Introduction} \label{Sec0}
The purpose of this paper is to initiate a systematic study of the irreducible
weight representations of direct limits of reductive
Lie algebras, and in particular of the classical simple direct limit Lie algebras
$A(\infty)$, $B(\infty)$, $C(\infty)$ and $D(\infty)$. 
We study  arbitrary, not 
necessarily highest weight, irreducible weight modules
and describe the supports of all such modules. 
The representation theory of the classical direct limit groups has been 
initiated in the pioneering works of G. Olshanskii \cite{O1} and \cite{O2}
and is now in an active phase, see the recent works of A. Habib, \cite{Ha}, 
L. Natarajan, \cite{Na}, K.-H. Neeb, \cite{Ne}, and 
L. Natarajan, E. Rodriguez-Carrington and J. A. Wolf, \cite{NRW}. 
Nevertheless, the structure theory of weight representations of the 
simple direct limit Lie algebras has until recently been still  in its infancy as 
only highest weight modules have been discussed in the literature, see the works of
Yu. A. Bahturin and G. Benkart, \cite{BB},  K.-H. Neeb, \cite{Ne}, and 
T. D. Palev, \cite{P}.

Our approach 
is based mainly on the recent paper \cite{DMP} in which a general method for
studying the support of weight representations of finite-dimensional 
Lie algebras (and Lie superalgebras) was developed. 
We prove first that the shadow of any irreducible weight module $M$ 
(over any root reductive direct limit Lie algebra, in particular
over $A(\infty)$, $B(\infty)$, $C(\infty)$
and $D(\infty)$) is well-defined, which means that for a 
given root $\alpha$ the intersection of the ray $\lambda + \R_+ \alpha$ with
the support of $M$, $\supp M$, is either finite for all $\lambda \in \supp M$
or infinite for all $\lambda \in \supp M$.  
Using this remarkable property of the support, we assign 
to $M$ a canonical parabolic subalgebra $\gp_M$ of $\gg$ and then compare $M$ 
with the irreducible quotient of a certain $\gg$-module induced from $\gp_M$. 
In the case of a finite-dimensional Lie algebra, the Fernando-Futorny parabolic
induction theorem, see \cite{Fe}, \cite{Fu} and \cite{DMP}, states that $M$ 
is always such a quotient and moreover that $\supp M$ is nothing but the 
support of the induced module. In the direct limit case we show that this is 
no longer true but nevertheless, using the fact that the shadow is 
well-defined, we obtain an explicit description of the support of any 
irreducible weight module. 

We discuss in more detail the following special cases: when $\supp M$ is finite
in all root directions (these are the finite integrable irreducible modules 
and they are analogues of finite-dimensional irreducible modules over 
finite-dimensional Lie algebras) and the more general case when $\supp M$ is 
finite in at least one of each two mutually opposite root directions. An 
interesting feature in the first case is that $M$ is not necessarily a highest
weight module, i.e. the analogues of finite-dimensional irreducible modules 
are already outside the class of highest weight modules. We present an 
explicit parametrization of all finite integrable irreducible modules. The 
modules corresponding to the second case are by definition pseudo highest 
weight modules and one of their interesting features is that in general they 
are not obtained by parabolic induction as in the Fernando-Futorny theorem.

Here is the table of contents:

1. Generalities on triangular decompositions and weight modules.

2. Direct limits of reductive Lie algebras.

3. The shadow of an irreducible weight module.

4. Integrable modules.

5. Pseudo highest weight modules.

6. The support of an arbitrary irreducible weight module.

{\bf Acknowledgement.} Discussions of various aspects of the subject of this
paper with Yu. A. Bahturin, G. Benkart, V. Futorny, Yu. I. Manin, O. Mathieu, 
G. I. Olshanskii, V. Serganova and J. A. Wolf have been very helpful to us. 
Special thanks are due to K.-H. Neeb who read an early version of the paper 
and pointed out some errors.
Both authors have been supported in part by an NSF GIG Grant, and I. Penkov 
acknowledges support from the Louis Pasteur University in Strasbourg where 
part of this work was done.

\section{Generalities on triangular decompositions and weight modules} 
\label{Sec1}
The ground field is $\C$. The signs $\cplus$ and
$\bcplus$ denote semi-direct sum of Lie algebras
(if $\gg = \gg' \cplus \gg''$, then $\gg'$ is an ideal in $\gg$ 
as well as a $\gg''$-module).
The superscript $^*$ always stands for dual space.
We set $\R_+ := \{r\in \R \, | \, r \geq 0 \}$, 
$\R_- := - \R_+$, $\Z_\pm := \Z \cap \R_\pm$, 
$\delta^{ij}$ is Kronecker's delta, and linear span with coefficients in
$\C$ (respectively, $\R$, $\Z$, $\R_\pm$) is denoted by
$< \quad >_\C$ (resp. by $< \quad >_\R$, $< \quad >_\Z$, 
$< \quad >_{\R_\pm}$). If $\gg$ is any Lie algebra and $M$ is a $\gg$-module,
we call $M$ {\it integrable} iff $\gg$ acts locally finitely on 
$M$; this terminology is introduced by V. Kac, in \cite{BB} integrable 
modules are called locally finite. 
If $\gg$ is a direct sum of Lie algebras, $\gg = \oplus_{s \in S} \gg^s$, 
and for each $s$ $M^s$ is an irreducible $\gg^s$-module with a
fixed non-zero vector $m^s \in M$, then $(\otimes_s M^s)(\otimes_s m^s)$ denotes 
the $\gg$-submodule of $\otimes_s  M^s$ generated by $\otimes_s m^s$.
It is easy to check that $(\otimes_s M^s)(\otimes_s m^s)$ is
an irreducible $\gg$-module.

Let $\gg$ be a Lie algebra. 
A {\it Cartan subalgebra of $\gg$} is by definition
a self-normalizing nilpotent Lie subalgebra $\gh \subset \gg$. 
We do not assume $\gg$ or $\gh$ to be finite-dimensional. 
A $\gg$-module $M$ is a {\it generalized weight $\gg$-module} iff
as an $\gh$-module $M$ 
decomposes as the direct sum $\oplus_{\lambda \in \gh^*} M^\lambda$, where
$$
M^\lambda := \{ \sv \in M \, | \, h-\lambda(h) {\rm {\, acts \, nilpotently \,
on \, }} \sv {\rm {\, for \, every \,}} h \in \gh\}.
$$
We call $M^\lambda$ the {\it generalized weight space of $M$ of weight 
$\lambda$}. Obviously, a generalized weight 
$\gg$-module is integrable as an $\gh$-module. 
The {\it support} of $M$, $\supp M$, consists 
of all weights $\lambda$ with $M^\lambda \neq 0$. A generalized weight module $M$ is
a {\it weight module} iff $\gh$ acts semi-simply on $M$, i.e. iff each generalized
weight space is isomorphic to the direct sum of one-dimensional $\gh$-modules.

Henceforth $\gg$ will denote a Lie superalgebra with a fixed proper Cartan 
subalgebra $\gh$ such that $\gg$ is a generalized weight module, i.e. 
\begin{equation} \label{eq1.0}
\gg = \gh \oplus (\oplus_{\alpha \in \gh^* \backslash \{0 \}} \gg^\alpha).
\end{equation}
The generalized weight spaces $\gg^\alpha$ are by definition the {\it 
root spaces} of $\gg$, $\Delta := \{ \alpha \in \gh^* \backslash
\{0\} | \gg^\alpha \neq 0 \}$ is the set of {\it roots of $\gg$}, and
the decomposition (\ref{eq1.0}) is the {\it root decomposition of $\gg$}.

%If $\gg$ is a direct sum of Lie
%algebras, $\gg = \oplus_{s \in S} \gg^s$, any given irreducible 
%generalized weight $\gg$-module $M$ can be described 
%in terms of irreducible generalized weight $\gg^s$-modules as follows. Fix
%an arbitrary non-zero generalized weight vector $m \in M$. For any 
%$s \in S$, set $\widehat{\gg^s}:=  \oplus_{t \neq s} \gg^t$ and let 
%$M^s$ and $\widehat{M^s}$ be respectively the $\gg^s$- and 
%$\widehat{\gg^s}$-submodules of $M$ generated by $m$. Using the fact 
%that $\gh^* = (\gh^s)^* \oplus (\widehat{\gh^s})^*$, where
%$\gh^s := \gh \cap \gg^s$ and $\widehat{\gh^s} := \gh \cap \widehat{\gg^s}$,  
%the reader will check that $M^s$ and $\widehat{M^s}$ 
%are irreducible modules. Furthermore, there is an obvious isomorphism
%of $\gg$-modules
%\begin{equation} \label{eq1.1}
%\varphi^s: M^s \otimes \widehat{M^s} \to M,
%\end{equation}
%such that $\varphi^s(m \otimes m)= m$. Consider now the infinite tensor
%product $\otimes_{s \in S} M^s$ with its coordinatewise 
%$\oplus_{s \in S} \gg^s$-module structure and let 
%$(\otimes_{s \in S}M^s)(m)$ be the 
%$\gg$-submodule of $\otimes_{s \in S} M^s$ generated by $\otimes_{s \in S}m$.
%$(\otimes_{s \in S} M^s)(m)$ is an irreducible  $\gg$-module because each
%of the modules $M^s$ is irreducible and every vector in \newline
%$(\otimes_{s \in S} M^s)(m)$ differs from $\otimes_{s \in S}m$ at finitely
%many coordinates only. Furthermore, it 
%is an exercise to check that $M \simeq (\otimes_{s \in S} M^s)(m)$. 

Below we recall the definitions of a triangular decomposition of $\Delta$
and of a Borel subalgebra for an arbitrary Lie algebra $\gg$ with a root
decomposition. For more details see \cite{DP2}.
A decomposition of $\Delta$
\begin{equation} \label{eq1.2}
\Delta = \Delta^+ \sqcup \Delta^-
\end{equation}
is a {\it triangular decomposition of } $\Delta$ iff the cone 
$<\Delta^+ \cup -\Delta^->_{\R_+}$ 
{\footnote {We use the term cone as a synonym for an 
$\R_+$-invariant additive subset of a real vector space.}} 
(or equivalently, its opposite
cone $<-\Delta^+ \cup \Delta^->_{\R_+}$) contains no (real) vector subspace.
Equivalently, (\ref{eq1.2}) is a triangular decomposition if the following is 
a  well-defined $\R$-linear partial order on $<\Delta>_\R$:
$$
\eta \geq \mu \Leftrightarrow \eta = \mu + \sum_{i=1}^n c_i \alpha_i
\text{ for some }  \alpha_i \in \Delta^+ \cup -\Delta^- \text { and } c_i 
\in \R_+, \text{ or } \mu = \eta .
$$
(A partial order is {\it $\R$-linear} if it is compatible with addition and
multiplication by positive real numbers, and if multiplication by negative
real numbers changes order direction.  In what follows, unless 
explicitly stated that it is partial, an order will
always be assumed linear, i.e. such that $\alpha \neq \beta$ implies
$\alpha < \beta$ or $\alpha > \beta$. In particular, an {\it $\R$-linear order}
is by definition an order which is in addition $\R$-linear.)
Any regular real hyperplane $H$ in $<\Delta>_\R$ (i.e. a codimension one 
linear subspace $H$ with $H \cap \Delta = \emptyset$)  
determines exactly two triangular decompositions of $\Delta$:
we first assign (in an arbitrary way) the sign $+$ to one of the two
connected components of $<\Delta>_\R \backslash H$, and the sign $-$ to
the other. $\Delta^\pm$ are then by definition the subsets of $\Delta$ which 
belong respectively to the "positive" and the "negative" connected components 
of $<\Delta>_\R \backslash H$. In general, not every triangular decomposition 
of $\Delta$ arises in this way, see \cite{DP1}. Nevertheless, it is true that 
every triangular decomposition is determined by a (not unique)
oriented maximal chain of vector suspaces in $<\Delta>_\R$: see the Appendix where we 
establish the precise interrelationship between 
oriented maximal chains in $<\Delta>_\R$, $\R$-linear orders on 
$<\Delta>_\R$ and triangular decompositions of $\Delta$.  

A Lie subalgebra $\bg$ of $\gg$ is by definition a {\it Borel 
subalgebra of} $\gg$ iff 
$\bg = \gh \bcplus ( \opp\limits_{\alpha \in \Delta^+} \gg^\alpha)$
for some triangular decomposition $\Delta = \Delta^+ \sqcup \Delta^-$.
In what follows a Borel subalgebra of $\gg$ always means a Borel subalgebra
containing the fixed Cartan subalgebra $\gh$.   
Adopting terminology from affine Kac-Moody algebras, we will call a Borel
subalgebra {\it standard} iff it corresponds to a triangular decomposition
which can be determined by a regular hyperplane $H$ in $<\Delta>_\R$.

If $\gb$ is a Borel subalgebra and $\lambda \in \gh^*$, the {\it Verma module
$\tilde{V}_\gb(\lambda)$ with $\gb$-highest weight $\lambda$} is by 
definition the induced module $U(\gg) \otimes_{U(\gb)} \C_\lambda$,
$\C_\lambda$ being the one-dimensional $\gb$-module on which $\gh$
acts via $\lambda$. Any quotient of $\tilde{V}_\gb(\lambda)$ is
by definition a {\it $\gb$-highest weight module}. Furthermore, 
$\tilde{V}_\gb(\lambda)$ has a unique maximal proper $\gg$-submodule
$I_\gb(\lambda)$, and $V_\gb(\lambda) := \tilde{V}_\gb(\lambda)/I_\gb(\lambda)$ 
is by definition the {\it 
irreducible $\gb$-highest weight $\gg$-module with highest weight $\lambda$}.

\section {Direct limits of reductive Lie algebras} \label{Sec2}
A homomorphism $\varphi: \gg \to \gg'$ of 
Lie algebras with root decomposition is a {\it root homomorphism}
iff $\varphi (\gh) \subset \gh'$ ($\gh$ and $\gh'$ 
being the corresponding fixed Cartan subalgebras) and $\varphi$ maps any 
root space of $\gg$ into a root space of $\gg'$. Let 
\begin{equation} \label{eq2.1}
\gg_1 @>{\varphi_1}>> \gg_2 @>{\varphi_2}>>  \ldots @>{\varphi_{n-1}}>> \gg_n
 @>{\varphi_n}>>  \ldots
\end{equation}
be a chain of homomorphisms of finite-dimensional Lie algebras and let 
$\gg := \dlim \gg_n$ be the direct limit Lie algebra. We say that 
$\gg$ is a {\it root direct limit of the system} (\ref{eq2.1}) iff 
$\varphi_n$ is a root homomorphism for every $n$. In the latter case $\gh_n$ 
denotes the fixed Cartan subalgebra of $\gg_n$. 
We define a Lie algebra $\gg$ to be a {\it root direct limit Lie algebra} 
iff $\gg$ is a root 
direct limit of some direct system  of the form (\ref{eq2.1}).
Furthermore, a non-zero Lie algebra 
$\gg$ is a {\it root reductive direct limit of the system} 
(\ref{eq2.1}) iff all $\gg_n$ are reductive, and $\gg$ is a {\it root simple 
direct limit Lie algebra} iff $\gg$ is a root direct limit of a 
system (\ref{eq2.1}) in which all $\gg_n$ are simple.

\begin{Pro} \label{Pro2.1} Let $\gg$ be a root direct limit Lie algebra. Then
$\gh := \dlim \gh_n$ is a Cartan subalgebra of $\gg$ and $\gg$ has a root 
decomposition with respect to $\gh$ such that $\Delta = \dlim \Delta_n$, where
$\Delta_n$ and $\Delta$ are respectively the roots of $\gg_n$ and $\gg$. 
If all root homomorphisms $\varphi_n$ are embeddings, one has simply 
$\gh = \cup_n \gh_n$ and $\Delta = \cup_n \Delta_n$.
\end{Pro}

{\it Proof.} A trivial exercise.
\qed

Every simple finite-dimensional Lie algebra $\gg$ 
is a root simple direct limit Lie algebra: we set $\gg_n := \gg$ and 
$\varphi_n := \id_{\gg}$. To define the simple Lie algebras $A(\infty)$, 
$B(\infty)$, $C(\infty)$ and $D(\infty)$ it suffices to let $\gg_n$ be the 
corresponding rank $n$ simple Lie algebra and to request that all $\varphi_n$ be
injective root homomorphisms, i.e. root embeddings. 
Indeed, one can check that in these cases the direct limit Lie algebra does 
not depend up to isomorphism on the choice of root embeddings $\varphi_n$. 
More generally, Theorem 4.4 in \cite{BB} implies that every 
infinite-dimensional root simple direct limit Lie
algebra is isomorphic to one of the Lie algebras $A(\infty)$, $B(\infty)$, 
$C(\infty)$ or $D(\infty)$. Note however, that for general root reductive direct limit Lie 
algebras the direct limit Lie algebra can 
depend on the choice of root
homomorphisms $\varphi_n$ even if they are embeddings.
Indeed, set for instance $\gg_n := A(2^n) \oplus B(2^n)$ and let
$\varphi_n, \varphi_n' : \gg_n \to \gg_{n+1}$ be root
embeddings such that $\varphi_n (A(2^n)) \subset A(2^{n+1})$
and $\varphi_n(B(2^n)) \subset B(2^{n+1})$ but 
$\varphi_n'(A(2^n) \oplus B(2^n)) \subset B(2^{n+1})$. Then 
$\gg \simeq A(\infty) \oplus B(\infty)$ but $\gg' \simeq B(\infty)$.

Let $\gg$ be a root reductive direct limit Lie algebra and let $\gk$ be a 
Lie subalgebra of $\gg$ such that $\gh \subset \gk$. 
Then $\Delta_\gk \subset \Delta_\gg$ and we can set 
$\Delta_\gk^\ss := \Delta_\gk \cap
(- \Delta_\gk)$. The Lie subalgebra $\gk^\ss := 
[\opp_{\alpha \in \Delta_\gk^\ss} \gk^\alpha,
\opp_{\alpha \in \Delta_\gk^\ss} \gk^\alpha]$ of $\gk$ is an
analogue of the semi-simple part of $\gk$ in the case when $\gk$ is not 
finite-dimensional. 

\begin{The} \label{The0}
Let $\gg$ be a root reductive direct limit Lie algebra. Then 

{\rm (i)} $\gg \simeq \gg^\ss \cplus A$, $A$ being an abelian Lie algebra of 
finite or countable dimension. Furthermore, $\gh =\gh^\ss \oplus A$, where
$\gh^\ss := \gh \cap \gg^\ss$. 

{\rm (ii)} $\gg^\ss \simeq \oplus_{s \in S} \gg^s$, $S$ being a finite or countable
family of root simple direct limit Lie algebras $\gg^s$.
\end{The}

{\it Proof.} (i) Decompose $\gh$ as $\gh = \gh^\ss \oplus A$ for some
vector space $A$. Then $A$ is an abelian subalgebra of $\gh$ 
(since $\gh$ itself is abelian) of at most countable dimension, and using 
Proposition 1 one checks that $\gg \simeq \gg^\ss \cplus A$.

(ii) Let $\gg = \dlim \gg_n$. Proposition \ref{Pro2.1} implies immediately that  
$\gg^\ss = \dlim \gg_n^\ss$.
Furthermore, since a Cartan subalgebra $\gh_n \subset
\gg_n$ is fixed for every $n$, there is a canonical decomposition 
as
$$
\gg_n = (\oplus_{t \in S_n} \gg^t) \oplus Z_n
$$ 
such that all $\gg^t$ are simple, $Z_n$ is abelian, and, for every 
$\alpha \in \Delta_n$ $\gg_n^\alpha \subset \gg^t$ for some
$t \in S_n$. Then, for any $t \in S_n$, either $\varphi_n(\gg^t) = 0$ or 
there is  $t' \in S_{n+1}$ so that $\varphi_n(\gg^t)$ is a non-trivial
subalgebra of $\gg^{t'}$
(we assume that the sets $S_n$ are pairwise disjoint). Put 
$S' = \cup_n \{ t \in S_n \,| \, \varphi_{n,n_1}(\gg^t) \neq 0 {\rm {\, for \,
every\,}} n_1 > n\}$, where $\varphi_{n_1, n_2}:= \varphi_{n_2} \circ \ldots \circ 
\varphi_{n_1}$ for $n_1 > n_2$ . Introduce an equivalence relation $\sim$ on 
$S'$ by setting  $t_1 \sim t_2$ for
$t_1 \in S_{n_1}$, $t_2 \in S_{n_2}$, iff there exists $n$ such that 
$\varphi_{n_1, n}(\gg^{t_1})$ and $\varphi_{n_2, n}(\gg^{t_2})$ belong to 
the same simple component of $\gg_n$. 
Define $S$ as the set of classes of $\sim$-equivalence.
For every $s \in S$ the set $\{\gg^t\}_{t \in s}$
is partially ordered by the maps 
$\varphi_{n_1,n_2}:\gg^{t_1} \to \gg^{t_2}$.
Let $\gg^s$ be the direct limit Lie algebra of a maximal 
chain of Lie algebras among $\{\gg^t\}_{t \in s}$ with respect to
this partial order.
Obviously $\gg^s$ is a root simple direct limit Lie algebra and it does not 
depend on the choice of the maximal chain. Finally, one checks easily that
$\gg^\ss = \dlim \gg_n^\ss \simeq \oplus_{s \in S} \gg^s$.
\qed

{\bf Example 1.} $gl(\infty)$ can be defined as the Lie algebra of all infinite
matrices $(a_{ij})_{i,j \in \Z_+}$ with finitely many non-zero entries. Then
$gl(\infty) \simeq A(\infty) \cplus \C$, but 
$gl(\infty) \not \simeq A(\infty) \oplus \C$ as the center of $gl(\infty)$ is trivial.

Theorem \ref{The0} gives an almost explicit description of all root reductive direct limit
Lie algebras. In particular, it implies that any root space of a root reductive 
direct limit Lie algebra has dimension one. 
Moreover, if $\pi : \gh^* \to (\gh^\ss)^*$ denotes the natural projection, then
Theorem \ref{The0} implies that $\pi$ induces a bijection between $\Delta$
and the set of roots of $\gg^\ss$. 
The exact relationship between irreducible generalized weight $\gg$-modules
and irreducible generalized weight $\gg^\ss$-modules (all of which turn
out to be automatically weight modules) is established in the following proposition.

\begin{Pro} \label{Pro2.2}
{\rm (i)} Every irreducible generalized weight $\gg$-module $M$ is a weight module.

{\rm (ii)} Every irreducible weight $\gg$-module $M$ is irreducible as a (weight) 
$\gg^\ss$-module. 

{\rm (iii)} Given any irreducible weight $\gg^\ss$-module $M^\ss$,
every $\lambda \in \gh^*$ with $\pi(\lambda) \in \supp M^\ss$ defines 
a unique structure of an irreducible weight $\gg$-module on $M^\ss$ which extends 
the $\gg^\ss$-module structure on $M^\ss$. If $M^\ss(\lambda)$ denotes 
the resulting $\gg$-module, then $M^\ss(\lambda) \simeq M^\ss(\lambda')$ iff
$\pi (\lambda - \lambda') = \sum_i c_i \pi(\alpha_i)$ with $\alpha_i \in \Delta$
implies $\lambda - \lambda' = \sum_i c_i \alpha_i$.
\end{Pro}

{\it Proof.} (i) Let $U^0$ denote the subalgebra of the enveloping algebra $U(\gg)$
generated by monomials of weight zero. Since $\gh$ acts semisimply on $\gg$, the 
symmetric algebra $S(\gh)$ belongs to the center of $U^0$. Furthermore, any generalized 
weight space $M^\lambda$ is an irreducible $U^0$-module and, by a general version of
Schur's Lemma, $S(\gh)$ acts via a scalar on $M^\lambda$. Therefore $M$ is a semi-simple 
$\gh$-module, i.e. $M$ is a weight module.

(ii) Follows immediately from (i) and Theorem \ref{The0}, (i).

(iii) Let any $a \in A$ act on the weight space $M^\mu$ of $M^\ss$ via
multiplication by $\lambda(a) + \sum_i \alpha_i(a)$, where 
$\mu - \pi(\lambda) = \sum_i c_i \pi(\alpha_i)$. It is straightforward to
verify that this equips $M^\ss$ with a well-defined $\gg$-module structure. 
The isomorphism criterion is also an easy exercise.
\qed

Our main objective in this paper is the study of the irreducible weight modules 
over root reductive direct limit Lie algebras. The case of finite-dimensional 
reductive Lie algebras
is discussed in particular in \cite{Fe}, \cite{Fu}, \cite{DMP}, \cite{M}. 

In the rest of this section we study the structure of the root simple
direct limit Lie algebras and $\gg$ stands for 
$A(\infty)$, $B(\infty)$, $C(\infty)$ or $D(\infty)$.
If $\vep_i$ are the usual linear functions on the Cartan subalgebras of
the simple finite-dimensional Lie algebras, see for example \cite{B} or 
\cite{Hu}, one can let $i$ run from $1$ to $\infty$ and then 
(using Proposition \ref{Pro2.1}) verify the following  list of roots:

$A(\infty) \quad : \quad \Delta = \{\vep_i-\vep_j \,|\, i\neq j\}$,

$B(\infty) \quad : \quad \Delta = \{\pm \vep_i,  \pm \vep_i \pm \vep_j \,|
\, i \neq j \}$,

$C(\infty) \quad : \quad \Delta = \{\pm 2\vep_i,  \pm \vep_i \pm \vep_j \,|
\, i \neq j \}$,

$D(\infty) \quad : \quad \Delta = \{\pm \vep_i \pm \vep_j \,|\, i \neq j \}$. \newline
Every sequence of complex numbers $\{\lambda^n\}_{n = 1, 2, \ldots}$ determines 
a weight $\lambda$ of $\gg$ by setting $\lambda(\vep_n) := \lambda^n$. 
For $\gg = B(\infty)$, $C(\infty)$ or
$D(\infty)$, every weight $\lambda$ of $\gg$ recovers the sequence
$\{\lambda^n := \lambda(\vep_n)\}_{n = 1, 2, \ldots}$;
for $\gg = A(\infty)$, the weight $\lambda$ recovers the sequence 
$\{\lambda^n := \lambda(\vep_n)\}_{n = 1, 2, \ldots}$ up to an additive constant only. 
Furthermore, $(\cdot, \cdot)$ denotes the bilinear form 
$( \cdot, \cdot): \gh^* \times \gh^* \to \C$ which extends the usual bilinear 
forms $(\cdot, \cdot)_n : \gh_n^* \times \gh_n^* \to \C$. If 
$\lambda \in \gh^*$ and $\alpha \in \Delta$, we set 
$\left< \lambda, \alpha\right>:= \frac{2(\lambda, \alpha)}{(\alpha,\alpha)}$. 
By definition, $\lambda \in \gh^*$ is an {\it integral weight} of $\gg$ iff 
$\left< \lambda, \alpha \right> \in \Z$ for every $\alpha \in \Delta$, and 
$\lambda$ is {\it $\gb$-dominant}, for a Borel subalgebra 
$\gb = \gh \oplus (\oplus_{\alpha \in \Delta^+} \gg^\alpha)$, 
iff $\left< \lambda, \alpha \right> \geq 0$ for all $\alpha \in \Delta^+$.
 
It is proved in \cite{DP2} (Proposition 2) that all Borel subalgebras of 
$A(\infty)$, $B(\infty)$, $C(\infty)$ an $D(\infty)$ are
standard. More precisely, for every Borel subalgebra $\gb$ of $\gg$ there is a 
linear function $\varphi : <\Delta>_\R 
\to \R$ such that  $\gb = \gh \oplus (\oplus_{\varphi(\alpha)>0} 
\gg^\alpha)$. Define an order (or a partial order) on the set 
$\{0\} \cup \{\pm \vep_i\}$ to be {\it $\Z_2$-linear} iff multiplication by $-1$
reverses the order. Then  Proposition 2 from \cite{DP2} is essentially equivalent
to the following statement.

\begin{Pro} \label{Pro2.3}
If $\gg = A(\infty)$, there is a bijection between Borel subalgebras of 
$\gg$ and orders on the set  $\{ \vep_i \}$.

If $\gg = B(\infty)$ or $C(\infty)$, there is a bijection between
Borel subalgebras of $\gg$ and $\Z_2$-linear orders on the set
$\{0\} \cup \{ \pm \vep_i \}$.

If $\gg = D(\infty)$, there is a bijection between  Borel subalgebras of 
$\gg$ and $\Z_2$-linear orders on the set $\{0\} \cup \{ \pm \vep_i \}$ 
with the property that if there is a minimal positive element with respect to
this order, this element is of the form $\vep_i$.
\end{Pro}

{\it Proof.} The pull-back via $\varphi$ of the standard order on $\R$
determines an order on $\{0\} \cup \Delta$ which induces
an order respectively on $\{\vep_i\}$ or $\{0\} \cup \{ \pm \vep_i \}$
as desired.
Conversely, for every order $\{\vep_i\}$, or respectively for every $\Z_2$-linear 
order on $\{0\} \cup \{\pm \vep_i\}$ as in the Proposition, there exists
a (non-unique) linear function $\varphi : <\{\vep_i\}>_\R \to \R$ such 
that $\Delta^+ = \varphi^{-1}(\R^+) \cap \Delta$.
\qed 

The result of Proposition 3 can be found in an equivalent form in \cite{Ne}. 
The case of $A(\infty)$ is due to V. Kac (unpublished).

{\bf Example 2.} $A(\infty)$ is naturally identified with the Lie algebra 
of traceless infinite matrices $(a_{i,j})_{i,j \in \Z_+}$ with
finitely many non-zero entries, as well as 
with the Lie algebra of traceless double infinite matrices  
$(a_{i,j})_{i,j \in \Z}$ with finitely many non-zero entries. 
The respective algebras of upper triangular matrices are non-isomorphic Borel 
subalgebras of $A(\infty)$ corresponding respectively to the orders 
$\vep_1 > \vep_2 > \vep_3 > \ldots$ and  
$\ldots > \vep_6 > \vep_4 > \vep_2 > \vep_1 > \vep_3 > \vep_5 > \ldots$.

A {\it parabolic subalgebra} of a root reductive direct limit Lie algebra is 
by definition a Lie subalgebra containing a Borel subalgebra. (The 
general definition of a parabolic subalgebra of a Lie algebra with
root decomposition is given in the Appendix.) Here is 
an explicit description of all parabolic subalgebras of $A(\infty)$,
$B(\infty)$, $C(\infty)$ and $D(\infty)$. 

\begin{Pro} \label{Pro2.4}
If $\gg = A(\infty)$, there is a bijection between parabolic subalgebras 
of $\gg$ and partial orders on the set  $\{ \vep_i \}$. 

If $\gg = B(\infty)$ or $C(\infty)$, there is a bijection between
parabolic subalgebras of $\gg$ and $\Z_2$-linear partial orders on the set
$\{0\} \cup \{ \pm \vep_i \}$.

If $\gg = D(\infty)$, there is a bijection between parabolic subalgebras of 
$\gg$ and $\Z_2$-linear partial orders on the set $\{0\} \cup \{ \pm \vep_i \}$ 
with the property that if $\vep_i$ is not comparable with $0$ for some $i$ 
(i.e. neither $\vep_i > 0$ nor $\vep_i <0$) then $\vep_j$ is also not
comparable with $0$ for some $j \neq i$.
\end{Pro}

{\it Proof.} 
Given a partial order on  $\{ \vep_i \}$ (respectively, a $\Z_2$-linear
partial
order on $\{0\} \cup \{ \pm \vep_i \}$ with the additional property for 
$\gg = D(\infty)$), it determines a unique partial order $>$ on the set
$\{0\} \cup \Delta$. Put 
$\gp_>:=\gh \oplus (\oplus_{\alpha > 0 \text{ or } \alpha \text{ not comparable with } 0} 
\gg^\alpha)$. Then  $\gp_>$ is the parabolic subalgebra corresponding 
to the initial partial order.
Conversely, let $\gp$ be a parabolic subalgebra. 
For $\alpha \in \Delta$, set $\alpha >_\gp 0$ iff $\gg^\alpha \subset \gp$ but
$\gg^{-\alpha} \not \subset \gp$. Using the explicit form of $\Delta$ it is easy
to verify that this determines a unique 
partial order on $\{ \vep_i \}$ (respectively on 
$\{0\} \cup \{ \pm \vep_i \}$ as desired).
\qed

The next proposition describes $\gp^\ss$ for any parabolic subalgebra $\gp$ of 
$\gg = A(\infty)$, $B(\infty)$, $C(\infty)$ or $D(\infty)$, 
and will be used in Section 3.

\begin{Pro} \label{Pro2.5} If $\gp \subset \gg$ is a parabolic subalgebra of $\gg$, 
then 

{\rm (i)} $\gp^\ss$ is isomorphic to a direct sum of simple Lie algebras each of
which is one of the following: \newline
{\rm -} $A(n)$ or $A(\infty)$, if $\gg= A(\infty)$; \newline
{\rm -} $A(n)$, $A(\infty)$, $B(n)$ or $B(\infty)$ with at most one simple
component isomorphic to $B(n)$ or $B(\infty)$, if $\gg = B(\infty)$; \newline
{\rm -} $A(n)$, $A(\infty)$, $C(n)$ or $C(\infty)$ with at most one simple
component isomorphic to $C(n)$ or $C(\infty)$, if $\gg = C(\infty)$; \newline
{\rm -} $A(n)$, $A(\infty)$, $D(n)$ or $D(\infty)$ with at most one simple
component isomorphic to $D(n)$ or $D(\infty)$, if $\gg = D(\infty)$.

{\rm (ii)} If $\gp^\ss \not \simeq 0$, then $\gp^\ss + \gh \simeq
(\oplus_{t \in T} \gk^t) \oplus Z$ , where $Z$ is abelian and $\gk^t$
is isomorphic to $gl(n)$ or $gl(\infty)$ for any $t \in T$ except at 
most one index $t_0 \in T$ for which $\gk^{t_0}$ is a root simple direct 
limit Lie algebra.
\end{Pro}

{\it Proof.} (i) Let, as in Proposition \ref{Pro2.4}, $>_\gp$ be the 
partial order corresponding to $\gp$. Partition the set $\{ \vep_i \}$ 
(respectively $\{0\} \cup \{ \pm \vep_i \}$) into subsets in such a way 
that two elements are comparable with respect to $>_\gp$ iff they belong to different subsets.
Denote by $S'$ the resulting set of subsets of $\{ \vep_i \}$ (respectively of
$\{0\} \cup \{ \pm \vep_i \}$). Let furthermore $S$ be a subset
of $S'$ which contains exactly one element of any pair $s, s'$ of mutually
opposite elements of $S'$, and for every $s \in S$ define
$\gg^s$ to be the Lie algebra generated by $\gg^\alpha$ where $\alpha \neq 0$
belongs to $s$ or is a sum of any two elements of $s$. 
It is straightforward to verify that $\gp^\ss \simeq \oplus_{s \in S, \gg^s \neq 0} \gg^s$ is 
the decomposition of $\gp^\ss$ into a direct sum of root simple direct limit
Lie algebras and that this decomposition satisfies (i).

(ii) The main point is to notice that if $s \neq -s$, then 
$\gg^s \simeq A(n)$ for some $n$ or $\gg^s \simeq A(\infty)$ and
furthermore, that there is at most one $s \in S$ such that $s = -s$.
To complete the proof it remains to show that each of the Lie
subalgebras $\gg^s$ for $s \neq -s$ can be extended to
a Lie subalgebra isomorphic to $gl(n)$ or $gl(\infty)$. 
This latter argument is purely combinatorial
and we leave it to the reader.
\qed

\section{The shadow of an irreducible weight module} \label{Sec3} 
Let $\gg$ be a Lie algebra with a root decomposition, 
$M$ be an irreducible generalized weight $\gg$-module and $\lambda$ be a fixed point
in $\supp M$. For any $\alpha \in \Delta$, consider the set $m_\alpha^\lambda
:= \{ q \in \R | \lambda + q \alpha \in \supp M \} \subset \R$. There are
four possible types of sets $m_\alpha^\lambda$: bounded in both directions; 
unbounded in both directions; bounded from above but unbounded from below; 
unbounded from above but bounded from below.
It is proved in \cite{DMP} that, if $\gg$ is finite-dimensional, 
the type of $m_\alpha^\lambda$ depends only on $\alpha$ and not on 
$\lambda$, and therefore the module $M$ itself determines
a partition of $\Delta$ into four mutually disjoint subsets:

$$
\Delta_M^F := \{ \alpha \in \Delta \, | \, m_\alpha^\lambda \, {\rm is \,
bounded \, in \, both \, directions} \},
$$
$$
\Delta_M^I := \{ \alpha \in \Delta \, | \, m_\alpha^\lambda \, {\rm is \,
unbounded \, in \, both \, directions} \},
$$
$$
\Delta_M^+ := \{ \alpha \in \Delta \, | \, m_\alpha^\lambda \, {\rm is \,
bounded \, from \, above \, and \, unbounded \, from \, below} \},
$$
$$
\Delta_M^- := \{ \alpha \in \Delta \, | \, m_\alpha^\lambda \, {\rm is \,
bounded \, from \, below \, and \, unbounded \, from \, above} \}.
$$
The corresponding decomposition 
\begin{equation} \label{eq3.1}
\gg = (\gg_M^F + \gg_M^I) \oplus \gg_M^+ \oplus \gg_M^-,
\end{equation}
where $\gg_M^F := \gh \oplus (\oplus_{\alpha \in \Delta_M^F} \gg^\alpha)$,
$\gg_M^I := \gh \oplus (\oplus_{\alpha \in \Delta_M^I} \gg^\alpha)$ and
$\gg_M^\pm := \oplus_{\alpha \in \Delta_M^\pm} \gg^\alpha$, is the
{\it $M$-decomposition} of $\gg$. The triple $(\gg_M^I, \gg_M^+, \gg_M^-)$
is the {\it shadow of $M$ onto $\gg$}.
If $\gg$ is infinite-dimensional we say that {\it the shadow of $M$ onto $\gg$
is well-defined} if it is true that the type of $m_\alpha^\lambda$ depends
only on $\alpha$ and not on $\lambda$ and thus the decomposition (\ref{eq3.1})
(as well as the triple $(\gg_M^I, \gg_M^+, \gg_M^-)$) is well-defined. 

In the case when $\gg$ is finite-dimensional and reductive, it is furthermore 
true that $\gp_M := (\gg_M^F + \gg_M^I) \oplus \gg_M^+$ is a parabolic 
subalgebra of $\gg$ whose reductive part is $\gg_M^F + \gg_M^I$, and that 
there is a natural surjection
$$
\varphi_M : U(\gg) \otimes_{U(\gp_M)} M^{\gg_M^+} \to M,
$$
where $M^{\gg_M^+}$ is the irreducible $(\gg_M^F + \gg_M^I)$-submodule
of $M$ which consists of all vectors in $M$ annihilated by  $\gg_M^+$.
Moreover, $\supp M$ simply coincides with \newline
$\supp (U(\gg) \otimes_{U(\gp_M)} M^{\gg_M^+})$. This is the Fernando-Futorny
parabolic induction theorem, see \cite {Fe} and also 
\cite{DMP}, and in particular it provides an 
explicit description of $\supp M$.

The main purpose of this paper is to understand analogues of these results 
for a root reductive direct limit Lie algebra $\gg$. Very roughly, the 
situation turns out to be as follows: the shadow of an arbitrary irreducible 
$\gg$-module $M$ exists and defines a parabolic subalgebra $\gp_M$ of $\gg$,
however the parabolic induction theorem does not hold. Nevertheless the 
existence of the shadow and the direct limit structure on $\gg$ enable us to 
obtain an explicit description of $\supp M$ for any $M$.

Throughout the rest of this paper (except the Appendix) $\gg:= \dlim \gg_n$ will be a root 
reductive direct limit Lie algebra and $M$ will be a fixed irreducible 
weight module over $\gg$.  
  
\begin{The} \label{The3.1}
The shadow of $M$ is well-defined.
\end{The}

{\it Proof.}
Note first that, for any given $\lambda \in \supp M$, there exist 
irreducible $(\gg_n + \gh)$-modules $M_n$ such that $\lambda \in \supp M_n$ for every $n$
and $\supp M = \cup_n \supp M_n$.
Indeed, fix $m \in M^\lambda$, $m \neq 0$, and define $M_n$ as any irreducible quotient of 
the $(\gg_n + \gh)$-module $U(\gg_n + \gh)\cdot m$. One checks 
immediately that $\supp M = \cup_n \supp M_n$.

Fix $\alpha \in \Delta$.
To prove the Theorem we need to show that if $m_\alpha^\lambda$ is bounded from 
above for some $\lambda \in \supp M$ (the proof for the case when $m_\alpha^\lambda$
is bounded from below is exactly the same) then $m_\alpha^\mu$ is bounded from above 
for any other $\mu \in \supp M$. First, let $\lambda' = \lambda + k \alpha$ be the end 
point of the $\alpha$-string through $\lambda$ in $\supp M$. Then fix $\mu \in \supp M$ 
and pick $N$ so that 
$\mu \in \supp M_N$, $\lambda' \in \supp M_N$ and $\alpha \in \Delta_N$. 
Consider now the $\alpha$-string through $\mu$. It is the union 
of the $\alpha$-strings through $\mu$ in $\supp M_n$ for $n = N, N+1, \ldots$.
The parabolic induction theorem implies that
$\supp M_n$ (and, in particular, the $\alpha$-string through $\mu$ in
$\supp M_n$) is contained entirely in an affine half-space in $\lambda +
<\Delta_n>_\R$ whose boundary (affine) hyperplane $H_n$ contains $\lambda'$ and 
is spanned by vectors in $\Delta_n$. Furthermore, $H_n' := H_n \cap 
(\lambda + <\Delta_N>_\R)$ is a hyperplane in $\lambda + <\Delta_N>_\R$
which contains $\lambda'$ and is spanned by vectors in $\Delta_N$.
Since the $\alpha$-string through $\mu$ belongs to $\lambda + <\Delta_N>_\R$,
the $\alpha$-string through $\mu$ in $\supp M_n$ is bounded by $H_n'$.
But there are only finitely many hyperplanes in $\lambda + <\Delta_N>_\R$
passing through $\lambda'$ and spanned by vectors in $\Delta_N$ and
hence among the hyperplanes $H_n'$ for $n = N, N+1, \ldots$ only finitely 
many are different. Therefore, the $\alpha$-strings through $\mu$ in $\supp M_n$ are
uniformly bounded from above, i.e. $m_\alpha^\mu$ is bounded from above.
\qed
  
The next theorem describes the structure of the $M$-decomposition and is 
an exact analogue of the corresponding theorem for finite-dimensional reductive
Lie algebras.

\begin{The} \label{The3.2}

{\rm (i)} 
$\gg_M^F$, $\gg_M^I$, $\gg_M^+$ and $\gg_M^-$ are Lie subalgebras of $\gg$ 
and  $\gh$ is a Cartan subalgebra for both $\gg_M^F$ and $\gg_M^I$.

{\rm (ii)}
$[(\gg_M^F)^\ss,(\gg_M^I)^\ss]=0$ and therefore 
$\gg_M^{FI}:= \gg_M^F + \gg_M^I$ is a Lie subalgebra of $\gg$.

{\rm (iii)}
$\gg_M^+$ and $\gg_M^-$ are $\gg_M^{FI}$-modules.

{\rm (iv)} 
$\gp_M := \gg_M^{FI} \oplus \gg_M^+$ and $\gg_M^{FI} \oplus \gg_M^-$ are 
(mutually opposite) parabolic subalgebras of $\gg$.
\end{The}

{\it Proof.} (i) Let $\alpha, \beta  \in \Delta$ such that  
$\alpha + \beta \in \Delta$. Lemma 2 in \cite{PS} implies that if 
$m_\alpha^\lambda$ and $m_\beta^\lambda$ are bounded from above, so is 
$m_{\alpha + \beta}^\lambda$. Furthermore, noting that 
$\supp M = \cup_n \supp M_n$ for some 
irreducible $(\gg_n+\gh)$-modules $M_n$ (see the proof of Theorem \ref{The3.1}) and 
that the support of every irreducible weight $(\gg_n +\gh)$-module 
is convex, we conclude that $\supp M$ is 
convex. Therefore if $m_\alpha^\lambda$ and $m_\beta^\lambda$ are 
unbounded from above, so is $m_{\alpha + \beta}^\lambda$. These two 
facts imply immediately that all four subspaces $\gg_M^F$, $\gg_M^I$, $\gg_M^+$ and $\gg_M^-$
are subalgebras of $\gg$. The fact that $\gh$ is a Cartan subalgebra for both $\gg_M^F$
and $\gg_M^I$ is obvious.

(ii) If $\alpha ' \in \Delta_M^F$ and $\beta ' \in \Delta_M^I$, then 
$\alpha ' + \beta ' \not \in \Delta$. Indeed, if $\alpha '+ \beta '\in \Delta$
and $m_{\alpha' + \beta'}^\lambda$ is bounded from above, then 
$m_{\beta '}^\lambda$ would be bounded from above because $\beta ' = 
-\alpha ' + (\alpha ' + \beta ')$.  If, on the other hand,
 $\alpha ' + \beta ' \in \Delta$ and $m_{\alpha ' + \beta '}^\lambda$ is 
unbounded from above, then $m_{\alpha '}^\lambda$ would be unbounded from 
above because $\alpha ' = -\beta ' + (\alpha ' + \beta ')$. Since both of these 
conclusions contradict the choice of $\alpha '$ and $\beta '$ we
obtain that $\alpha ' +\beta ' \not \in \Delta$ and thus that  
$[(\gg_M^F)^\ss,(\gg_M^I)^\ss]=0$.

(iii) If $\alpha ' \in \Delta_M^F$, $\beta ' \in \Delta_M^+$
and $\alpha ' +\beta ' \in \Delta$, then again $m_{\alpha ' + \beta '}^\lambda$
is bounded from above. Assuming that $m_{\alpha ' + \beta '}^\lambda$
is bounded from below, we obtain that $m_{\beta '}^\lambda$ is bounded
from below as well since $-\beta' = -(\alpha ' + \beta ') + \alpha '$, 
which contradicts the fact that $\beta ' \in \Delta_M^+$. Hence 
$m_{\alpha ' + \beta '}^\lambda$ is unbounded from below and $\alpha ' +\beta '
\in \Delta_M^+$. This proves that $\gg_M^+$ is a $\gg_M^F$-module. One shows 
in a similar way that $\gg_M^\pm$ are $\gg_M^F$- and $\gg_M^I$-modules.

(iv) is a direct corollary of (i), (ii) and (iii).
\qed

We define an irreducible weight $\gg$-module $M$ to be
{\it cuspidal} iff $\gg = \gg_M^I$.

In the rest of this section we prove that for every parabolic subalgebra $\gp$
of $\gg$ there is an irreducible $\gg$-module $M$ so that 
$\gp = \gp_M$. Indeed, there is the following more general

\begin{The} \label{The3.3} Let $\gg$ be a root reductive direct limit Lie 
algebra. For any given splitting $\Delta = \Delta^F \sqcup \Delta^I \sqcup \Delta^+
\sqcup \Delta^-$ with $\Delta^F = -\Delta^F$, $\Delta^I = -\Delta^I$,
$\Delta^- = -\Delta^+$, and such that its corresponding decomposition 
\begin{equation} \label{eq3.2}
\gg = (\gg^F +\gg^I) \oplus \gg^+ \oplus \gg^-
\end{equation}
satisfies the properties {\rm (i) - (iii)} of Theorem {\rm \ref{The3.2}}, 
there exists an irreducible weight module 
$M$ for which {\rm (\ref {eq3.2})} is the $M$-decomposition of $\gg$. 
\end{The}

{\it Proof.} 
We start with the observation that it suffices to prove the Theorem for
a root simple direct limit Lie algebra. Indeed, let as in Theorem \ref{The0}
$\gg \simeq (\oplus_{s \in S} \gg^s) \cplus A$ and assume that for each
$s$ $M^s$ is an irreducible weight $\gg^s$-module corresponding as in the Theorem 
to the restriction of the decomposition (\ref{eq3.2}) to $\gg^s$. Fix 
non-zero vectors $m^s \in M^s$.
Then the reader will verify straightforwardly that, for any pair 
$(M^\ss := (\otimes_s M^s)(\otimes_s m^s), \lambda)$ as in Proposition \ref{Pro2.2}, 
(iii), the $\gg$-module $M:= M^\ss(\lambda)$ is as required by the Theorem.
Therefore in the rest of the proof we will assume that $\gg$ is a root simple 
direct limit Lie algebra. We will prove first that cuspidal modules exist.
 
\begin{Lem} \label{Lem3.1}
Let $\gg$ be a simple finite-dimensional Lie algebra. 
For any given weight $\mu \in \gh^*$, 
there exists a cuspidal $\gg$-module $M^I$ such that the center $Z$ of $U(\gg)$
acts on $M$ via the central character $\chi_\mu : Z \rightarrow \C$ (obtained 
by extending $\mu$ to a homomorphism $\tilde \mu : S^\cdot (\gh) \rightarrow 
\C$ and composing with Harish-Chandra's homomorphism $Z\rightarrow 
S^\cdot(\gh)$).
\end{Lem}

{\it Proof of Lemma \ref{Lem3.1}.} It goes by induction on the rank of $\gg$.
It is a classical fact that $sl(2)$ admits a cuspidal module of any given 
central character, so it remains to make the induction step.
Let $\gg'$ be a reductive subalgebra of $\gg$ which contains $\gh$ and such
that $\rk \gg' = \rk \gg - 1$. Let $M'$ be a cuspidal $\gg'$-module with 
central character $\chi_{\mu'}'$, where $\mu' \in \gh^*$ and $w'(\mu') - \mu \not\in
<\Delta>_\Z$ for any $w'$ in the Weyl group $W'$ of $\gg'$. 
Denote by $U'$ the subalgebra of $U(\gg)$ generated by $U(\gg')$
and $Z$. Since $U(\gg)$ is a free $Z$-module (Kostant's Theorem) and 
$U(\gg') \cap Z = \C$, $U'$ is isomorphic to $U(\gg')\otimes_\C Z$.
Therefore the tensor product $M_{\chi_\mu}' := M' \otimes \C_{\chi_\mu}$,  
$\C_{\chi_\mu}$ being the 
1-dimensional $Z$-module corresponding to $\chi_\mu$, is a well-defined $U'$-module.
Consider now any 
irreducible quotient $M$ of the induced $\gg$-module $U(\gg) \otimes _{U'}
M_{\chi_\mu}'$. Obviously $M$ is a weight module and we claim that $M$ is 
cuspidal. Assuming the contrary, $M$ would be a quotient of 
$U(\gg) \otimes_{U(\gp)} M''$ where $\gp \supset \gg'$ and $M''$ is a 
cuspidal $\gg'$-module of central character $\chi_\mu'$. Then, since $U(\gg)$ is an 
integrable $U(\gg')$-module and $M'$ would be a subquotient of 
$U(\gg) \otimes_{U(\gp)} M''$, we would have $w'(\mu') - \mu \in <\Delta>_\Z$ for some
$w' \in W'$, which 
is contradiction. Therefore $M$ is cuspidal.
\qed

If $\gg = \gg^I$, $\gg = \dlim \gg_n$ being an infinite-dimensional root simple  
direct limit Lie algebra, we can assume that 
$\rk \gg_{n+1} - \rk \gg_n = 1$.
We then construct $M^I$ as the direct limit $\dlim M_n$,
each $M_n$ being a cuspidal $\gg_n$-module build by induction precisely as 
in the Proof of Lemma \ref{Lem3.1}. This proves the Theorem in the cuspidal case.

Let now $\gg \neq \gg^I$. According to Theorem \ref{The0}, 
$\gg^F \simeq \oplus_{s \in S} \gg^s \cplus A^F$ and 
$\gg^I \simeq \oplus_{t \in T} \gg^t \cplus A^I$, where each of the algebras 
$\gg^r$ for  $r \in S \cup T$ is a root simple direct limit Lie algebra and 
$A^F$ and $A^I$ are respectively abelian subalgebras of $\gg^F$ and $\gg^I$.  
Since the Theorem is proved for the cuspidal case, 
we can choose a cuspidal irreducible $\gg^t$ module $M^t$ for each $t \in T$
and fix a non-zero vector $m^t \in M^t$. Then $\widehat{M^I} :=
(\otimes_{t \in T} M^t) (\otimes_t m^t)$ is a cuspidal $(\gg^I)^\ss$-module. 
Let $M'$ denote $\widehat{M^I}$ considered as a $(\gg^F + \gg^I)^\ss$-module
with trivial action of $(\gg^F)^\ss$.

If $\dim \gg < \infty$, we can assume furthermore, by Lemma \ref{Lem3.1},
that the $(\gg^I)^\ss$-module $\widehat{M^I}$ is chosen in such a way 
that the central character of $\widehat{M^I}$ corresponds to
an orbit $W^I \cdot \eta$ ($W^I$ being the Weyl group of $(\gg^I)^\ss$)
which contains no weight of the form $\eta + \kappa$ for an integral weight
$\kappa$ of $(\gg^I)^\ss$. Then it is not difficult to verify 
(we leave this to the reader) that this condition ensures that if $M^{FI}$ 
denotes $\widehat{M^I}$ with its obvious $(\gg^F + \gg^I) \oplus \gg^+$-module
structure and if $M$ is the irreducible quotient of 
$U(\gg) \otimes_{U((\gg^F + \gg^I) \oplus \gg^+)} M^{FI}$, the $M$-decomposition of $\gg$
is nothing but (\ref{eq3.2}).
Theorem \ref{The3.3} is therefore proved for a finite-dimensional reductive
$\gg$. 

A direct generalization of this argument does not go through for
an infinite-dimensional $\gg$. Instead, there is the following lemma
which allows us to avoid referring to central characters.

\begin{Lem} \label{Lem3.2}
If $\gg$ is infinite-dimensional, 
the $(\gg^F + \gg^I)^\ss$-module structure on $M'$ can be extended to a 
$(\gg^F +\gg^I)$-module structure in such a 
way that $\left< \mu', \alpha \right> \notin \Z$ for some (and hence any) 
$\mu' \in \supp M'$ and any $\alpha \in \Delta^+ \sqcup \Delta^-$. 
\end{Lem}
{\it Proof of Lemma \ref{Lem3.2}.} 
Theorem \ref{The3.3}, (iv) implies that $\gg_M^F + \gg_M^I = (\gp_M)^\ss + \gh$.
We will present the proof in the case when $\gg_M^F + \gg_M^I \neq \gh$.
The case when $\gg_M^F + \gg_M^I \neq \gh$ is dealt with in a similar way.

Using Proposition \ref{Pro2.5}, (ii) we conclude that
$\gg_M^F + \gg_M^I \neq \gh$ implies
\begin{equation} \label{eq3.10} 
\gg_M^F + \gg_M^I \simeq (\oplus_{r \in R} \gg^r) \oplus \gg^0 \oplus Z, 
\end{equation}
where each $\gg^r$ is isomorphic to $gl(n)$ for some $n$ or to $gl(\infty)$, 
$R$ is a finite or countable set which does not contain $0$ (and may be empty),
$\gg^0$ is a root simple direct limit Lie algebra, and $Z$ is abelian. 

Let $\Delta^s$ be the root system of $\gg^s$ for $s \in R \cup \{0\}$. Set
$^s\Z := \{ k \in (\Z_+ \backslash \{0\}) | {\text { there is }} \newline l \in 
(\Z_+ \backslash \{0\}) {\text { with }} 
\vep_k - \vep_l \in \Delta^s \}$ and  $^\#\Z := (\Z_+ \backslash \{0\}) 
\backslash (\cup_{s \in R \cup \{0\}} \, ^s\Z)$. 
As explained in section \ref{Sec1}, any sequence 
$\{\lambda^n\}_{n=1, 2, \ldots} \subset \C$ determines a weight $\lambda$
of $(\gg^F + \gg^I)^\ss$, however $\lambda$ reconstructs 
only the subsequence
$\{\lambda^n\}_{n \in \cup_{r \in R \cup \{0\}}\, ^r\Z}$ 
up to an additive constant for every index from $R$.
Fix now $\mu \in \supp M'$ and let $\{\mu^n\}_{n=1, 2, \ldots}$
be a sequence which determines $\mu$. To prove the Lemma, it suffices to find
constants $c^i \in \C$ for $i \in \Z$ with the following properties: \newline
- $c^k = c^l$ whenever $k$ and $l$ belong to one and the same set $^s\Z$; \newline
- $c^k = 0$ for $k \in \, ^0\Z$; \newline 
- $\mu^k + c^k \notin \Z$ for every $k \in (\Z_+ \backslash \{0\})$ for which 
$\vep_k  \in \Delta^+ \sqcup \Delta^-$; $2 (\mu^k + c^k) \notin \Z$ for every 
$k \in (\Z_+ \backslash \{0\})$ for which $2 \vep_k  \in \Delta^+ \sqcup \Delta^-$;
and $\pm (\mu^k + c^k) \pm (\mu^l + c^l) \notin \Z$ for every 
$k, l \in (\Z_+ \backslash \{0\})$ for which 
$\pm \vep_k \pm \vep_l \in \Delta^+ \sqcup \Delta^-$.
\newline 
Then $\mu'$ will be the weight of $(\gg^F +\gg^I)$ determined by the
sequence $\{\mu^n + c^n\}_{n=1, \infty}$.

Consider the set of sequences $\{ c^k \}_{k \in P \subset (\Z_+ \backslash \{0\})}$  
with the three properties as above. 
Introduce an order $<$ on this set by putting
$\{ {c'}^{k'} \}_{k' \in P'} < \{ {c''}^{k''}\}_{k'' \in P''}$ iff 
$P' \subset P''$ and ${c'}^k = {c''}^k$ for every $k \in P'$. The reader will
check that every chain (with respect to $<$) of sequences is bounded and that for every  
$\{ c^k \}_{k \in P}$ with $P \neq (\Z_+ \backslash \{0\})$ there is a sequence greater than 
$\{ c^k \}_{k \in P}$. 
Therefore any maximal element, which exists by Zorn's lemma, is a sequence $\{c^k\}$ 
with the required properties.
Lemma \ref{Lem3.2} is proved.
\qed

Assuming that $\gg$ is infinite-dimensional,
 let now $M^{FI}$ be  $M'$ considered as a $(\gg^F + \gg^I) \oplus \gg^+$-module 
with trivial action of $\gg^+$. Then we define $M$ as the (unique) irreducible 
quotient of the $\gg$-module $U(\gg) \otimes_{U((\gg^F + \gg^I) \oplus \gg^+)}M^{FI}$. 
Using Lemma \ref{Lem3.2}, the reader will verify that the $M$-decomposition of 
$\gg$ is precisely the decomposition (\ref{eq3.2}).
The proof of Theorem \ref{The3.2} is therefore complete.
\qed

\section{Integrable modules} \label{Sec4}
\begin{Pro} \label{Pro4.1} Let $M$ be integrable. Then

{\rm (i)} $\gg = \gg_M^F + \gg_M^I$; furthermore, $\gg \neq \gg_M^F$ implies 
$\dim M^\lambda = \infty$ for any $\lambda \in \supp M$.

{\rm (ii)} If $\gg$ is infinite-dimensional and simple,  $M$ is
either cuspidal or $\gg = \gg_M^F$, and both cases are possible.
\end{Pro}

{\it Proof.} (i) As a consequence of the integrability of $M$, $\supp M$ is 
$W$-invariant. Thus $\Delta_M^\pm = \emptyset$ and $\gg = \gg_M^F + \gg_M^I$.
Furthermore, for any $\alpha \in \Delta$ the subalgebra $\gg^{\R \alpha}$ of $\gg$ generated
by $\gg^\alpha$ and $\gg^{-\alpha}$ is isomorphic to $sl(2)$ and the integrability
of $M$ implies that as a $\gg^{\R \alpha}$- module $M$ is isomorphic to a direct sum of finite-dimensional
modules. The assumption that $\dim M^\lambda < \infty$ for some $\lambda \in \supp M$ would
lead us to the conclusion that $\lambda$ belongs to the support of only finitely many 
of the finite-dimensional $\gg^{\R \alpha}$-modules, which would mean that $\supp M$ is finite 
in the direction of both $\alpha$ and $-\alpha$. Therefore $\gg \neq \gg_M^F$ implies
$\dim M^\lambda = \infty$ for any $\lambda \in \supp M$.

(ii) If both $\gg_M^F \neq \gh$ and $\gg_M^I \neq \gh$, Theorems \ref{The0} 
and \ref{The3.2} imply that $(\gg_M^F)^\ss$ and $(\gg_M^I)^\ss$ are proper
ideals in $\gg$, which is a contradiction. Thus $\gg_M^F = \gh$ or
$\gg_M^I = \gh$, i.e. $\gg = \gg_M^I$ or $\gg= \gg_M^F$. 

Clearly, $\gg = \gg_M^F$ if $M = \gg$ is the adjoint representation,
so the case $\gg = \gg_M^F$ is obviously possible. 
To prove that integrable cuspidal  modules exist, it is enough to
construct a tower of irreducible finite-dimensional $\gg_n$-modules $M_n$,
$\ldots \to M_n \to M_{n+1} \to \ldots$, such that
the support of $M_n$ is shorter than the support of $M_{n+1}$ in all 
root directions of $\gg_n$. Obviously, $M:= \dlim M_n$ is then an integrable 
cuspidal $\gg$-module.
Here is an explicit example for $\gg = A(\infty)$.
Let $\lambda_n = n\vep_1 - n\vep_2$ and let $\gb_n$ be the Borel subalgebra
of $\gg_n$ corresponding to the order $\vep_2 < \vep_3 < \ldots <\vep_n <\vep_1$.
Set $M_n := V_{\gb_n}(\lambda_n)$ for $n \geq 2$. It is obvious that the support
of $M_{n+1}$ is longer that the support of $M_n$ in the direction of every root 
of $\gg_n$. The reader will also verify that $M_n$ appears in the 
decomposition of $M_{n+1}$ as a direct sum of $\gg_n$-modules and hence
an embedding $M_n \to M_{n+1}$ exists.
The proof of Proposition \ref{Pro4.1} is therefore complete.
\qed

Obviously, if $M$ is integrable and is a highest weight module with respect
to some Borel subalgebra $\gb$, then $\gg = \gg_M^F$.
The existence of cuspidal integrable modules  is a significant 
difference with the case of a finite-dimensional Lie 
algebra, where any integrable irreducible module is finite-dimensional and is thus a 
highest weight module with respect to every Borel subalgebra. In contrast, for
$\gg = A(\infty)$, $B(\infty)$, $C(\infty)$ and $D(\infty)$, 
the trivial $\gg$-module is the only 
irreducible $\gg$-module which is a highest weight module for all Borel
subalgebras. Moreover, as we will see in Example 3 below,
the equality $\gg = \gg_M^F$ does not guarantee the 
existence of a Borel subalgebra with respect to which $M$ is a highest weight 
module.
We define $M$ to be {\it finite integrable} iff $\gg = \gg_M^F$. The rest of this section
is devoted to the study of finite integrable irreducible $\gg$-modules $M$
over an arbitrary root reductive direct limit Lie algebra $\gg$.

The simplest type of finite integrable modules are highest weight integrable modules
and they are studied in the recent papers \cite{BB}, \cite{NRW} and \cite{Ne}. 
In \cite{BB} the highest weights of integrable highest weight modules are 
computed explicitly.
In \cite{NRW} 
the integrable modules $V_\gb(\lambda)$ appear in the context of Borel-Weil-Bott's 
theorem and are realized as the unique non-zero  cohomology groups of
line bundles on $G/B$. In \cite{Ne} irreducible highest weight modules with respect
to general Borel subalgebras
are considered and in particular it is proved that their integrability is equivalent to 
unitarizability.
In Theorem \ref{The4.20} below we discuss highest weight modules and in particular 
establish a direct limit version of H. Weyl's character formula for integrable highest 
weight modules.

First we need to recall the notion of  basis of a Borel subalgebra of $\gg$.
A subset $\Sigma \subset \Delta^+$ is a {\it basis} of 
$\bg = \gh \oplus (\oplus_{\alpha \in \Delta^+} \gg^\alpha)$  
iff $\Sigma$ is a linearly independent set and
every element of $\Delta^+$ is a linear combination of elements of $\Sigma$ with
non-negative integer coefficients; the elements of $\Sigma$ are then the {\it simple
roots of $\gb$}. 
Not every Borel subalgebra admits a basis. For root simple direct limit
Lie algebras a basis of $\gb$ is the same as a weak basis in the terminology
of \cite{DP2}, and in \cite{DP2} all Borel subalgebras admitting a weak basis 
are described. The result of \cite{DP2} implies that $\gb$ admits a basis iff
the corresponding order on $\{\vep_i\}$, or respectively on $\{0\} \cup
\{\pm \vep_i\}$, has the following property: for every pair of elements of
$\{\vep_i\}$ (respectively of $\{0\} \cup \{\pm \vep_i\}$) there are only finitely 
many elements between them. This latter criterion has been established 
also by K.-H. Neeb in \cite{Ne}.

\begin{The} \label{The4.20} Let $M = V_\gb (\lambda)$. 

{\rm (i)} $M$ is integrable iff $\lambda$
is an integral $\gb$-dominant weight. Furthermore, 
if $M$ is integrable, then $\supp M = C^\lambda$, where 
$C^\lambda$ is the 
intersection of the convex hull of $W \cdot \lambda$ with 
$\lambda +  <\Delta>_\R$.

{\rm (ii)} $M \simeq V_{\gb'}(\lambda')$ for given $\gb'$ and  $\lambda'$ iff 
there exists $w \in W$ for which $\lambda' = w(\lambda)$ and there is a parabolic subalgebra
$\gp$ of $\gg$ containing both $w(\gb)$ and $\gb'$ such that the $\gp$-submodule of 
$M$ generated by $M^{\lambda'}$ is one-dimensional.

{\rm (iii)} $\dim M^\mu < \infty$ for all $\mu \in \supp M$ iff 
$M \simeq V_{\tilde{\gb}}(\lambda)$ for some Borel subalgebra $\tilde{\gb}$ of $\gg$
which admits a basis.

{\rm (iv)} If $M$ is integrable, $\gb$ admits a basis, and $\ch M:= \sum_{\mu \in \supp M}
\dim M^\mu \cdot e^\mu$ is the formal character of $M$, we have
\begin{equation} \label{eqW}
D \, \cdot \, \ch M = \sum_{w \in W} (\sgn w) e^{w(\lambda + \rho_\gb) - \rho_\gb},
\end{equation}
where $D= \Pi_{\alpha \in \Delta^+} (1 - e^{-\alpha})$ and $\rho_\gb \in \gh^*$ is a weight
for which $\rho_\gb(\alpha) = 1$ for all simple roots of $\gb$.
\end{The}

{\it Proof.} (i) An exercise. In \cite{BB} a more general criterion for the integrability of
$V_\gb(\lambda)$ is established
for simple direct limit Lie algebras which are not necessarily root simple 
direct limit Lie algebras.

(ii) Define $\gp_\lambda \supset \gb$ as the maximal parabolic subalgebra of 
$\gg$ such that the $\gp_\lambda$-submodule of $M$ generated by $M^\lambda$ is 
one-dimensional. Then  
$M \simeq V_{\gb''}(\lambda)$ iff $\gb''$ is
a subalgebra of $\gp_\lambda$. Finally, it is an exercise to check that
$M \simeq V_{\gb'}(\lambda')$ iff
$\lambda' = w(\lambda)$ and $\gb'$ is a subalgebra of $\gp_{\lambda'} = w(\gp_\lambda)$
for some $w \in W$.

(iii) If $\tilde{\gb}$ admits a basis, then, using the Poincare-Birkhoff-Witt theorem, 
one verifies that all weight spaces of the Verma module $\tilde{V}_\gb(\lambda)$ are 
finite-dimensional (for any $\lambda \in \gh^*$), and hence all weight spaces 
of $V_\gb(\lambda)$ are finite-dimensional as well. 

Conversely, let $M= V_\bg(\lambda)$
be a highest weight module with finite-dimensional weight spaces.
We need to prove the existence of $\tilde{\gb}$ which admits a basis
and such that $M \simeq V_{\tilde{\gb}}(\lambda)$.
Let $\gp_\lambda$ be as above and 
let $\gg \simeq (\oplus_{s \in S} \gg^s) \cplus A$ as in Theorem \ref{The0}.
Clearly, there is a Borel subalgebra $\tilde{\gb} \subset \gp_\lambda$ 
which admits a basis iff for every $s \in S$ there is a 
Borel subalgebra $\tilde{\gb^s}$ of $\gg^s$ admitting basis such that
$\tilde{\gb^s} \subset \gp_\lambda^s := \gp_\lambda \cap \gg^s$.
Furthermore, $M^s := V_{\gb^s} (\lambda |_{\gh^s})$ (where
$\gb^s := \gb \cap \gg^s$ and $\gh^s := \gh \cap \gg^s$) is an irreducible
$\gg^s$-submodule of $M$. Thus it suffices to prove (iii) for the root simple
direct limit Lie algebras $\gg^s$ and their highest weight modules $M^s$.
For a finite-dimensional $\gg^s$ (iii) is trivial, so 
we need to consider only the case when $\gg^s = A(\infty)$,
$B(\infty)$, $C(\infty)$, $D(\infty)$. 

Let $>_{\gp_\lambda^s}$ be the partial order corresponding to $\gp_\lambda^s$ and
let $\alpha \in (\Delta^s)^+$ ($\Delta^s$ being the root system of $\gg^s$) be a 
difference of two elements $\delta_1$ and $\delta_2$ of $\{\vep_i\}$ 
(respectively of $\{0\} \cup \{\pm \vep_i\}$). 
Then it is not difficult to check that, if there are infinitely many elements of 
$\{\vep_i\}$ (respectively of $\{0\} \cup \{\pm \vep_i\}$) between
$\delta_1$ and $\delta_2$ with respect to $>_{\gp_\lambda^s}$, the weight space of 
$M^s$ with weight $\lambda - \alpha$ is infinite-dimensional.
Therefore, the assumption that all Borel subalgebras $\tilde{\gb}^s$,
such that $M^s \simeq V_{\tilde{\gb}^s}$ admit no basis is contradictory,
i.e. a Borel subalgebra $\tilde{\gb}^s$ exists as required.

(iv) To prove formula (\ref{eqW}) it suffices to notice that if $\gb$ admits
a basis, $\gg$ can
be represented as $\gg = \dlim \gg_n$, where each of the simple roots of
$\gb_n = \gb \cap \gg_n$ is a simple root of $\gb_{n+1} = \gb \cap \gg_{n+1}$. 
Then  $M = \dlim V_{\gb_n} 
(\lambda|_{\gh_n})$ and in both sides of (\ref{eqW}) terms of the form
$c \cdot e^{\lambda + \mu}$ for $\mu \in <\Delta_n>_\R$ 
appear only as they appear in the respective sides of (\ref{eqW})
 for the $\gg_n$-module $V_{\bg \cap \gg_n} (\lambda|_{\gh_n})$.  
Therefore, Weyl's original formula implies the infinite version (\ref{eqW}). 
\qed

We now turn our attention to general finite integrable modules $M$.
Here is an example of a finite integrable $M$ 
which is not a highest weight module with respect to any Borel
subalgebra of $\gg$.

{\bf Example 3.} Let $\gg = A(\infty)$ and $\bg \subset \gg$ be the 
Borel subalgebra of $\gg$ corresponding to the order 
$\vep_1 > \vep_2 > \vep_3 > \ldots$. Set $\lambda_n := \vep_1 + \cdots +
\vep_n - n \vep_{n+1}$. Since  
$\lambda_n \in \supp V_{\gb_{n+1}}(\lambda_{n+1})$ and the weight space  
$V_{\gb_{n+1}}(\lambda_{n+1})^{\lambda_n}$ is one-dimensional, 
there is a unique (up to a multiplicative constant) embedding of weight 
$\gg_n$-modules 
$V_{\gb_n}(\lambda_n) \to V_{\gb_{n+1}}(\lambda_{n+1})$. Set 
$M := \dlim V_{\bg_n}(\lambda_n)$. Then $\gg = \gg_M^F$ and all weight spaces 
of $M$ are one-dimensional but $M$ is not a highest weight module with respect 
to any Borel subalgebra of $\gg$.

The following two theorems provide an explicit parametrization of all finite integrable
modules as well as an explicit description of their supports. Let $W_n \subset W$
denote the Weyl group of $\gg_n$.

\begin{The} \label{The4.1} Let $M$ be finite integrable. 

{\rm (i)} $M \simeq \dlim M_n$ for some direct system of finite-dimensional
irreducible $\gg_n$-modules $M_n$.

{\rm (ii)} $\supp M$ determines $M$ up to isomorphism.

{\rm (iii)} Fix a Borel subalgebra $\gb$ of $\gg$. Then  
$M \simeq \dlim V_{\bg_n}(\lambda_n |_{\gh_n})$ for some sequence
$\{ \lambda_n \}$ of integral weights of $\gg$ such that $\lambda_n|_{\gh_n}$
is a $\gb_n$-dominant weight of $\gg_n$ and 
$\lambda_n$ belongs to an edge of the convex hull of $W_{n+1} \cdot \lambda_{n+1}$.
\end{The}

{\it Proof.} (i) Fix $\lambda \in \supp M$.
We claim that there is a unique irreducible finite-dimensional $\gg_n$-module 
$M_n^\lambda$ such that 
$\supp M_n^\lambda = (\lambda + <\Delta_n>_\Z) \cap \supp M$. Indeed,
the equality $\gg = \gg_M^F$ implies that 
$(\lambda + <\Delta_n>_\Z) \cap \supp M$ is a finite set and hence there is 
$n'$ such that $(\lambda + <\Delta_n>_\Z) \cap \supp M$ is contained in 
the support of an irreducible finite-dimensional $\gg_{n'}$-module $M_{n'}'$
constructed as in the proof of Theorem \ref{The3.1}. Moreover, the
construction of $M_{n'}'$ implies that $(\lambda + <\Delta_n>_\Z)
\cap \supp M_{n'}' = (\lambda + <\Delta_n>_\Z) \cap \supp M$. But, as
the reader will check, there is a unique irreducible finite-dimensional
$\gg_n$-module $M_n^\lambda$ with 
$\supp M_n^\lambda = (\lambda + <\Delta_n>_\Z) \cap \supp M_{n'}'$.
Finally $M_n^\lambda$ is a $\gg_n$-submodule of $M$ according to its
construction, and furthermore $M = \dlim M_n^\lambda$.
  
(ii) The crucial point is that (according to its construction) $M_n^\lambda$ 
depends only on $\supp M$ and $\lambda \in \supp M$. Furthermore, for any
$\lambda' \in \supp M$, there is $m$ such that $\lambda' \in \supp M_m^\lambda$
and hence $M_n^\lambda \simeq M_n^{\lambda'}$ for $n>m$. Letting now $n$ go to
$\infty$, and noting that for any pair $\lambda , \lambda' \in \supp M$ there
is a compatible system of isomorphisms $M_n^\lambda \simeq M_n^{\lambda'}$
for all $n>m$, we conclude that 
$M = \dlim M_n^\lambda \simeq \dlim M_n^{\lambda'}$, i.e. that $\supp M$
determines $M$ up to isomorphism.

(iii) Since the module $M_n^\lambda$ defined in (i) is finite-dimensional,
$M_n^\lambda \simeq V_{\gb_n} (\widetilde{\lambda_n})$ for some 
$\widetilde{\lambda_n} \in \gh_n^*$.
There is a unique weight $\mu$ of $\gg$ such that $\mu \in <\Delta_n>_\R$
and $\mu|_{\gh_n} = \lambda|_{\gh_n} - \widetilde{\lambda_n}$. 
Set $\lambda_n := \lambda - \mu \in \gh^*$.
Then $M \simeq \dlim M_n^\lambda \simeq V_{\gb_n} (\widetilde{\lambda_n})
= \dlim V_{\gb_n}(\lambda_n|_{\gh_n})$ and
$\lambda_n$ belongs to an edge of the convex hull of $W_{n+1} \cdot \lambda_{n+1}$
because $\supp M_n^\lambda = (\lambda + <\Delta_n>_\Z) \cap \supp M$.
\qed

\begin{The} \label{The4.2}
Fix a Borel subalgebra $\bg$ of $\gg$. Let $\{ \lambda_n \} \subset \gh^*$ 
be a sequence of integral weights of $\gg$ such that $\lambda_n |_{\gb_n}$ is
a $\gb_n$-dominant weight of $\gg_n$ and 
$\lambda_n$ belongs to an edge of the convex hull of $W_{n+1} \cdot \lambda_{n+1}$.
Then 

{\rm (i)} $M$ is finite integrable.

{\rm (ii)} $\supp M = \cup_n C^{\lambda_n}$.

{\rm (iii)} $M$ is a highest weight module with respect to some Borel 
subalgebra of $\gg$ iff there is $n_0$ so that 
$\lambda_n \in W \cdot \lambda_{n_0}$ for any $n \geq n_0$.

{\rm (iv)} If $M' \simeq \dlim V_{\gb_n}(\lambda_n'|_{\gh_n})$, 
$\{\lambda_n'\}$ being
a sequence of integral weights of $\gg$, such that $\lambda_n' |_{\gh_n}$ is
a $\gb_n$-dominant weight of $\gg_n$ and 
$\lambda_n'$ belongs to an edge of the convex hull of $W_{n+1} \cdot \lambda_{n+1}'$
$M \simeq M'$ iff
there is $n_0$ so that $\lambda_n = \lambda_n'$ for $n \geq n_0$.
\end{The}
{\it Proof.} The proof is not difficult and is left to the reader.
\qed

{\bf Example 4.}
Let $\gg = A(\infty)$ and $M = \gg$ be the adjoint module. If 
$\gb$ is the Borel subalgebra corresponding to the order $\vep_1 > \vep_2 > \ldots$, 
then $M \simeq \dlim V_{\bg_n}(\vep_1 - \vep_n)$. Since $W \cdot (\vep_1 - \vep_n) = 
\Delta$, Theorem \ref{The4.2}, (iii) implies that there exists a Borel subalgebra
$\gb$ such that the adjoint module is a 
$\gb$-highest weight module. 
Furthermore, it is not difficult to verify that all such Borel subalgebras $\gb$ are
precisely the Borel subalgebras which correspond to orders on $\{ \vep_i\}$
for which there exists a pair of indices $i_0, j_0$ so that $\vep_{i_0} > \vep_i$
and $\vep_i > \vep_{j_0}$ for all $i \neq i_0$, $i \neq j_0$.

\section{Pseudo highest weight modules} \label{Sec5}
In the case of a finite-dimensional Lie algebra, an irreducible weight module 
$M$ with $\gg_M^I = \gh$ is necessarily a highest weight module for some
Borel subalgebra, see \cite{DMP}. As we already know (Example 3)
this is no longer true for the direct limit algebras we consider.
We define a {\it pseudo highest weight module} as an irreducible weight module $M$
such that $\gg_M^I = \gh$. Pseudo highest weight modules provide counterexamples also
to the obvious extension of the parabolic induction theorem to root 
reductive direct limit Lie algebras. 
(The module $M$ from Example 3 above does not provide such a counterexample since
in this case $\gp_M = \gg$ and $M = U(\gg) \otimes_{U(\gp_M)} M$.) 
Indeed, it suffices to construct $M$ with $\gg_M^F = \gg_M^I = \gh$ such
that $M$ is not a highest weight module with respect to $\gh \oplus \gg_M^+$
and therefore admits no surjection of the form 
$U(\gg) \otimes_{U(\gp_M)} M' \to M$ for an irreducible $\gp_M$-module
$M'$. Here is such an example.

{\bf Example 5.} Let $\gg = A(\infty)$, $B(\infty)$, $C(\infty)$ or $D(\infty)$.
Fix a Borel subalgebra $\gb \subset \gg$ and 
let $\lambda \in \gh^*$ be such that $\left< \lambda, \alpha 
\right> \not \in \Z$ for all $\alpha \in \Delta$. Construct $\{ \lambda_n \}$
inductively by setting $\lambda_3 := \lambda$, $\lambda_{n+1} := \lambda_n 
+ \alpha$, where $\alpha_n$ is a simple root of $\gb_{n+1}$ 
which is not a root of $\gb_n$.
There are obvious embeddings $V_{\gb_n}(\lambda_n|_{\gh_n}) \to  
V_{\gb_{n+1}}(\lambda_{n+1}|_{\gh_{n+1}})$ so we can set
$M := \dlim V_{\gb_n}(\lambda_n|_{\gh_n})$. Then 
$\gg = \gh \oplus \gg_M^+ \oplus \gg_M^-$ where $\gh \oplus \gg_M^+ = \gb$. 
Therefore the only Borel subalgebra with respect to which $M$ could be 
a highest weight module is $\gb$. But a direct verification shows that $M$
has no non-zero $\gb$-highest weight vector and is therefore not a 
$\gb$-highest weight module.

For any irreducible weight module $M$ there is the following natural question: 
what are the integrable root directions of $M$, i.e. for which
roots $\alpha \in \Delta$ is  $M$ an integrable $\gg^\alpha$-module? 
It is a remarkable fact that if $\Delta_M^\int$ is the set of all roots $\alpha$ 
for which $M$ is $\gg^\alpha$-integrable, then $\gp_M^\int := \gh \oplus 
(\oplus_{\alpha \in \Delta_M^\int} \gg^\alpha)$ is a Lie subalgebra of $\gg$, 
and $\gp_M^\int$ is nothing but the subset of elements in $\gg$ which 
act locally finitely on $M$.
For a finite-dimensional Lie algebra this can be proved using Gabber's theorem, 
\cite{G}, see Corollary 2.7 in \cite{Fe} or Proposition 1 in \cite{PS}. For
a root direct limit Lie algebra the statement 
follows immediately from the case of a finite-dimensional Lie algebra.

\begin{Pro} \label{Pro5.11}

{\rm (i)} $\gp_M^\int = (\gg_M^F + (\gg_M^I \cap \gp_M^\int)) \oplus \gg_M^+$.

{\rm (ii)} If $M$ is a pseudo highest weight module, then $\gp_M^\int = \gp_M$;
in particular $\gp_M^\int $ is a parabolic subalgebra of $\gg$.
\end{Pro}

{\it Proof.}
(i) It is obvious that $\gg_M^F \subset \gp_M^\int$ and $\gg_M^+ \subset \gp_M^\int$.
Furthermore, $\gp_M^\int \cap \gg_M^- = 0$. Indeed, assuming that 
$\gg_M^- \cap \gp_M^\int \neq 0$ we would have that $\gg^\alpha$ acts locally
nilpotently on $M$ for some $\alpha \in \Delta_M^-$, which is impossible
as then $\supp M$ would have to be invariant with respect to the reflection 
$\sigma_\alpha \in W$. Since $\gp_M^\int \supset \gh$, $\gp_M^\int$ is a weight 
submodule of $\gg$ and $\gp_M^\int = \gp_M^\int \cap ((\gg_M^F + \gg_M^I) \oplus \gg_M^+)
= (\gg_M^F +(\gg_M^I \cap \gp_M^\int)) \oplus \gg_M^+$.

(ii) The statement follows immediately from (i) since $\gg_M^I = \gh \subset \gg_M^F$ for
a pseudo highest weight module $M$.
\qed

The following Proposition provides a more explicit description of $\gp_M = \gp_M^\int$
for highest weight modules. It is a version of the main result of \cite{DP2} 
adapted to highest weight modules with respect to arbitrary Borel subalgebras
a root reductive Lie algebra $\gg$.
If $M = V_\gb(\lambda)$, we say that  a root $\alpha \in \Delta^-$ is {\it $M$-simple} iff 
$\alpha = \beta + \gamma$ with $\beta, \gamma \in \Delta^-$ 
implies $\left<\lambda, \beta\right> = 0$ or
$\left<\lambda, \gamma \right> = 0$.

\begin{Pro} \label{Pro5.12}
Let $M=V_\gb(\lambda)$ and let $\Sigma_{\lambda,\gb}^\int$ be the set of all 
$M$-simple roots $\delta \in \Delta^-$ for which $M$ is 
$\gg^\delta$-integrable. Then, for any $\alpha \in \Delta^-$,
$M$ is $\gg^\alpha$-integrable (equivalently, $\alpha \in \Delta^- \cap \Delta_M^F$) iff 
$\alpha \in <\Sigma_{\lambda, \gb}^\int>_{\Z_+}$.
\end{Pro}

{\it Proof.} 
Proposition \ref{Pro5.11}, (ii) implies that if $\alpha \in 
<\Sigma_{\lambda, \gb}^\int>_{\Z_+}$ then  $M$ is $\gg^\alpha$-integrable.
Let, conversely,  $M$ be $\gg^\alpha$-integrable.  
If $\alpha$ is not $M$-simple, then $\alpha = \beta + \gamma$ for some 
$\beta, \gamma \in \Delta^-$ with 
$\left< \lambda, \beta \right> \neq 0$ and $\left< \lambda, \gamma \right> \neq 0$.
The reader will then check that $|\left< \lambda, \alpha \right>|$ is strictly 
bigger than both 
$|\left< \lambda, \beta \right>|$ and $|\left< \lambda, \gamma \right>|$.
Choose $n$ big enough so that $\alpha, \beta ,\gamma \in \Delta_n$.
Applying the main Theorem from \cite{DP2} to $V_{\gb_n}(\lambda|_{\gh_n})$,
we obtain that $V_{\gb_n}(\lambda|_{\gh_n})$ is both
$\gg^{\beta}$-integrable and $\gg^{\gamma}$-integrable. Therefore $M$ is also
$\gg^\beta$-integrable
as well as $\gg^\gamma$-integrable.
To conclude that $\alpha \in <\Sigma_{\lambda, \gb}^\int>_{\Z_+}$ whenever 
$V_\gb(\lambda)$ is $\gg^{\alpha}$-integrable one applies now induction on 
$|\left< \lambda, \alpha \right>|$.
\qed

Although the parabolic induction theorem does not hold for
pseudo highest weight modules, their supports can be described explicitly.
In the next section we prove a general theorem describing the support of
any irreducible weight module. An open question about pseudo highest
weight modules is whether the statement (ii) of Theorem \ref{The4.1}
extends to any pseudo highest weight module, i.e. whether such a module is 
determined up to an isomorphism by its support.

\section{The support of an arbitrary irreducible weight module} \label{Sec 6}
Let $M$ be an arbitrary irreducible weight module $M$ with corresponding partition 
$\Delta = \Delta_M^F \sqcup \Delta_M^I \sqcup \Delta_M^+ \sqcup \Delta_M^-$.
Define the {\it small Weyl group $W^F$ of} $M$ as the Weyl group of $\gg_M^F$.  
For  $\lambda \in \gh^*$, set $K_M^\lambda := (W^F \cdot \lambda) + 
<\Delta_M^I>_\Z + <\Delta_M^->_{\Z_+}$ and $K_M^{\lambda, n} := 
((W^F \cap W_n) \cdot \lambda) + <\Delta_M^I \cap \Delta_n>_\Z + <\Delta_M^- 
\cap \Delta_n>_{\Z_+}$.

\begin{Lem} \label {Lem6.1}
For any  $\lambda_0 \in \supp M$ and any $n$, there exists $\lambda \in \supp M$ 
such that \newline $(\lambda_0 + <\Delta_n>_\R) \cap \supp M = K_M^{\lambda, n}$.
\end{Lem}

{\it Proof.} If $\Delta_n \subset \Delta_M^I$, set $\lambda = \lambda_0$.
Assume now that  $\Delta_n \not \subset \Delta_M^I$.
Consider the cone $K:=<\Delta_M^F \cup
\Delta_M^+>_{\R_+}$. Then $(\lambda_0 + K) \cap (\supp M \cap <\Delta_n>_\R)$ 
is a finite set. As in the proof of Theorem \ref{The3.1},
$\supp M = \cup_N \supp M_N$ for some irreducible $\gg_N$-modules $M_N$. 
Furthermore, there is
$N_0$ for which $\supp M_{N_0} \supset (\lambda_0 + K) \cap (\supp M
\cap <\Delta_n>_\R)$. Let $\lambda \in \supp M_{N_0} \cap <\Delta_n>_\R$ be such that
$\lambda + \alpha \not \in \supp M_{N_0}$ for any $\alpha \in \Delta_n \cap (\Delta_M^F 
\cup \Delta_M^+)$. (Such a weight $\lambda$ exists because otherwise we would have $\Delta_n \subset
\Delta_M^I$.) Clearly then $(\lambda_0 + <\Delta_n>_\R) \cap \supp M = 
K_M^{\lambda, n}$.
\qed

Fix now a Borel subalgebra $\gb \subset (\gg_M^F +\gg_M^I) \oplus 
\gg_M^+$ of $\gg$ and let $\lambda \in \supp M$. 
Then Lemma \ref{Lem6.1} enables us to construct a sequence
$\lambda_3, \lambda_4, \ldots, \lambda_n, \ldots$ such that $\lambda_3 = \lambda$
and $\lambda_n$ is such that \newline
- $\lambda_n + \alpha \not \in \supp M$ for any $\alpha \in 
(\Delta_M^F \cup \Delta_M^+) \cap \Delta^+$; \newline
- $(\lambda_1 + <\Delta_n>_\R) \cap \supp M = K_M^{\lambda, n}$, for 
$n \geq 4$. \newline
This sequence reconstructs $\supp M$ and is, in a certain sense,
assigned naturally to $M$. For a precise
formulation, define an equivalence relation $\sim_M$ on $\gh^*$ by setting 
$\lambda_1 \sim_M \lambda_2$ iff there exists $w \in W^F$ such that 
$w(\lambda_1) - \lambda_2 \in <\Delta_M^I>_\Z$. 
Let furthermore $\gh_M^*$ denote the set of $\sim_M$-equivalence classes 
and let $p: \gh^* \mapsto \gh_M^*$ be the natural projection.

The following theorem is our most general result describing $\supp M$
for an arbitrary irreducible weight $\gg$-module $M$. It is a straightforward
corollary of the construction of the sequence $\{\lambda_n\}$.

\begin{The} \label{The6.1}
{\rm (i)} $K^{\lambda_1} \subset K^{\lambda_2} \subset \ldots $ and 
$\supp M = \cup_n K_M^{\lambda_n}$.

{\rm (ii)} The sequence $\{p(\lambda_n)\}$ depends on $\lambda$ only, and
if $\lambda'$ is another element of $\supp M$, there is $n_0$ so that
$p(\lambda_n) = p(\lambda_n')$ for all $n > n_0$.

{\rm (ii)} If $\hat{M}$ is an irreducible weight $\gg$-module with 
corresponding sequence $\{\widehat{\lambda_n}\}$, 
then $\supp M = \supp \hat M$ iff the there is $n_0$ such that $p(\lambda_n) = 
p(\widehat{\lambda_n})$ for $n > n_0$.
\end{The}
\qed

The existence of a sequence of weights $\{\lambda_n\}$ for which 
$\supp M = \cup_n K_M^{\lambda_n}$ follows directly from the fact
that $\supp M = \cup_n \supp M_n$ for some irreducible $\gg_n$-modules $M_n$
(see the proof of Theorem \ref{The3.1}). The main advantage of constructing
the sequence $\{\lambda_n\}$ is that $M$ stably determines the sequence 
$\{p(\lambda_n)\}$. It is an open question to describe explicitly all possible
sequences $\{p(\lambda_n)\}$, or equivalently all possible supports of 
irreducible weight $\gg$-modules with a given shadow.  

\section* {Appendix. Borel subalgebras of $\gg$ and chains of subspaces in $<\Delta>_\R$} 
\label{SecAp}
In this Appendix $\gg$ is an arbitrary Lie algebra with a fixed Cartan 
subalgebra $\gh$ such that $\gg$ admits a root decomposition (\ref{eq1.0}).
Here we establish a precise interrelationship 
between Borel subalgebras $\gb$ of $\gg$ containing $\gh$ (which we call
simply Borel subalgebras), $\R$-linear orders on $<\Delta>_\R$,
and oriented maximal chains of vector subspaces in $<\Delta>_\R$.
Our motivation is the following. If
$\gg$ has finitely many roots (i.e. for instance, if $\gg$ is finite-dimensional),
every Borel subalgebra is defined by a (non-unique) regular hyperplane $H$ in 
$<\Delta>_\R$, i.e. by a hyperplane $H$ such that $H \cap \Delta = \emptyset$. 
If $\gg$ has infinitely many roots this is known to be no longer true (i.e.
in general there are Borel subalgebras which do not correspond to any regular 
hyperplane in $<\Delta>_\R$ {\footnote {This is the case for any Kac-Moody algebra
which is not finite-dimensional.}}) and in \cite{DP1} we have shown that in the case
when $<\Delta>_\R$ is finite-dimensional every Borel subalgebra can be defined 
by a maximal flag of linear subspaces in $<\Delta>_\R$. However, when $<\Delta>_\R$
is infinite-dimensional, the situation is more complicated and deserves a careful formulation.

We start by stating the relationship between Borel subalgebras $\gb$ with 
$\gb \supset \gh$ and $\R$-linear orders on $<\Delta>_\R$. 

\begin{Pro} \label{ProAp1}
Every $\R$-linear order on $<\Delta>_\R$ determines a unique Borel subalgebra, 
and conversely, every Borel subalgebra is determined by an (in general not unique) 
$\R$-linear order on $<\Delta>_\R$.
\end{Pro}

{\it Proof.} If $>$ is an $\R$-linear order on $<\Delta>_\R$, set $\Delta^\pm :=
\{ \alpha \in \Delta \, | \, \pm \alpha >0\}$. Then $\Delta^+ \sqcup \Delta^-$ 
is a triangular decomposition and thus $>$ determines the Borel subalgebra
$\gh \oplus (\oplus_{\alpha \in \Delta^+} \gg^\alpha)$. The proof of the converse statement is
left to the reader (it is a standard application of Zorn's Lemma).
\qed

In what follows we establish a bijection between $\R$-linear orders on 
$<\Delta>_\R$ and oriented maximal chains of vector subspaces in $<\Delta>_\R$.
If $V$ is any real  vector space, a  {\it chain of vector subspaces} of $V$ is defined as a 
set of subspaces $F = \{F_\alpha\}_{\alpha \in A}$ of $V$ such that $\alpha \neq \beta$ 
implies a proper inclusion $F_\alpha \subset F_\beta$ or $F_\beta \subset F_\alpha$.
The set $A$ is then automatically ordered.
A chain $F$ is a {\it flag} iff as an ordered set $A$ can be identified with
a (finite or infinite) interval in $\Z$.
A chain $F$ is {\it maximal} iff it is not properly contained in any other
chain. Maximal chains of linear subspaces may be somewhat counterintuitive as the 
following example shows that a linear space of countable dimension admits non-countable 
maximal chains and vise versa.

{\bf Example 6.} Let $V$ be a countable-dimensional space 
with a basis $\{e_r\}_{r \in \Q}$. The chain of all subspaces
$V_t := <\{e_r \, | \, r < t\}>_\R$ for $t \in \R$,
$V_r' := <\{e_r \, | \, r \leq r\}>_\R$ for $r \in \Q$, $V_{-\infty}:=0$ and
$V_\infty := V$ is a maximal chain
of cardinality continuum.
Let $U := \C[[x]]$ be the space of formal power series in the indeterminate
$x$. The dimension of $U$ is continuum,
 however $\{0\} \cup F$, $F$ being the flag  $\{F_i := x^i W\}_{i \in \Z_+}$, 
provides an example of a countable maximal
chain in $U$.

Furthermore, it turns out that every maximal chain is determined uniquely by
a certain subchain. Define a chain $F$ to be {\it basic} 
iff it is minimal with the following property: for every $x \in V$ there exists a pair
$F_\alpha \subset F_\beta \in F$ with $\dim F_\beta/F_\alpha =1$
such that $x \in F_\beta \backslash F_\alpha$.

\begin{Lem} \label{LemA1}
Every maximal chain in $V$ admits a unique basic subchain and conversely,
every basic chain is contained in a unique maximal chain.
\end{Lem}

{\it Proof.} Let $F = \{F_\alpha\}_{\alpha \in A}$ be a maximal chain in $V$.
For every non-zero $x \in V$ set $F_x := \cup_{x \not \in F_\alpha} F_\alpha$ and
$F_x' := F_x \oplus \R x$. Then $F \cup \{F_x, F_x'\}$ is a chain in $V$
and hence $F_x , F_x' \in F$. Let $G := \cup_{0 \neq x \in V} \{F_x , F_x'\}$.
Then $G$ is obviously a basic chain contained
in $F$. Noting that $G$ consists exactly of all pairs of subspaces from $F$ 
with relative codimension one, we conclude that $G$ is the unique basic chain contained in $F$.

Let now $G = \{G_\beta\}_{\beta \in B}$ be a basic chain in $V$ and 
let $<$ be the corresponding order on $B$. 
A subset $C$ of $B$ is a {\it cut} of $B$ if
$\alpha \in C$ and $\beta < \alpha$ implies $\beta \in C$, and if
$\alpha \in C$ and $\beta \not \in C$ implies $\beta > \alpha$.
The set $\bar{B}$ of all cuts of $B$ is naturally ordered.
Set $H_\alpha := \cup_{\beta \in \alpha} G_\beta$ for $\alpha \in \bar{B}$.
One checks immediately that $H = \{H_\alpha\}_{\alpha \in \bar{B}}$
is a maximal chain in $V$ which contains $G$. On the other hand, given any 
maximal chain $F = \{F_\alpha\}_{\alpha \in A}$ in $V$
containing $G$, one notices that $F_\alpha = \cup_{G_\beta \subset F_\alpha} G_\beta$.
Then the map $\varphi: A \to \bar{B}$, $\varphi(\alpha) := \{\beta \, | \,
G_\beta \subset F_\alpha\}$ is an embedding of $F$ into $H$ and hence $F = H$.
The Lemma is proved.
\qed

An {\it orientation of a basic chain} $G$ is a labeling of 
the two half-spaces of $G_\beta \backslash G_\alpha$
for every pair of indices $\alpha < \beta$ with $\dim G_\beta / G_\alpha =1$,
by mutually opposite signs $\pm$. An {\it orientation of a maximal chain} $F$
is an orientation of its basic subchain. A maximal chain $F$ is {\it oriented}
iff an orientation of $F$ is fixed.

\begin{The} \label{TheAp}
There is a bijection between oriented maximal chains in 
$<\Delta>_\R$ and $\R$-linear orders on $<\Delta>_\R$.
\end{The}

{\it Proof.} Lemma \ref{LemA1} implies that it is enough to establish
a bijection between oriented basic chains in $<\Delta>_\R$ and 
$\R$-linear orders on $<\Delta>_\R$.

Let $G$ be an oriented basic chain in $<\Delta>_\R$.
We define an $\R$-linear order $>_G$ on $<\Delta>_\R$ by
setting $x >_G 0$ or $x <_G 0$ according to the sign of the  
half-space of $G_\beta \backslash G_\alpha$ to which $x$ belongs.
Conversely, given an $\R$-linear order $>$ on $<\Delta>_\R$ we build the corresponding 
oriented basic chain as follows. For a non-zero vector $x \in V$
we set $V_x := \{ y \in V \, | \, \pm(c y + x) > 0 {\rm {\, iff \,}} \pm x > 0 \}$. 
Then for every pair of non-zero vectors  $ x , y \in V$ exactly one of the 
following is true: $V_x = V_y$, $V_x \subset V_y$ or $V_x \supset V_y$. 
Therefore the set of distinct subspaces among $\{ V_x \}$ is a chain $F$. 
The very construction of $F$ shows that it is a basic chain in $<\Delta>_\R$.
To complete the proof of the Theorem it remains to orient this chain in the obvious way.
\qed

Finally we define a decomposition of $\Delta$,
\begin{equation} \label{eqAp10}
\Delta = \Delta^- \sqcup \Delta^0 \sqcup \Delta^+,
\end{equation}
to be {\it parabolic} iff $\Pi(\Delta^-) \cap \Pi(\Delta^+) = \emptyset$,
$0 \not \in \Pi(\Delta^\pm)$ and 
$\Pi(\Delta) \backslash \{0\} = \Pi(\Delta^-) \sqcup
\Pi(\Delta^+)$ is a triangular decomposition of $\Pi(\Delta) \backslash \{0\}$
(i.e. the cone $<\Pi(\Delta^+) \cup - \Pi(\Delta^-)>_{\R_+}$ contains no
vector subspace),
where $\Pi$ is the projection $<\Delta>_\R \to <\Delta>_\R / <\Delta^0>_\R$.
If $\Delta^0 = \emptyset$, a parabolic decomposition is a triangular decomposition. 
Given a parabolic decomposition (\ref{eqAp10}), its corresponding {\it parabolic
subalgebra} is by definition $\gp:= \gh \oplus (\oplus_{\alpha \in \Delta^0 \sqcup
\Delta^+} \gg^\alpha)$.

Given a linear subspace $V'$ in $V$, we define a {\it $V'$-maximal} chain of vector 
subspaces in $V$ to be the preimage in $V$ of a maximal chain in $V/V'$. 
An $V'$-maximal chain is {\it oriented} iff the corresponding maximal chain in 
$V/V'$ is oriented. Generalizing the corresponding construction for Borel subalgebras,
we have

\begin{The} \label{theAp10}
Any parabolic decomposition {\rm (\ref{eqAp10})} is determined by
some oriented $<\Delta^0>_\R$-maximal chain of vector subspaces in $V$, 
and conversely any oriented $V'$-maximal chain in $<\Delta>_\R$, for
an arbitrary subspace $V'$ of $<\Delta>_\R$, defines a unique parabolic decomposition of 
$\Delta$ with $\Delta^0 = \Delta \cap V'$.
\end{The}

{\it Proof.} An exercise. 
\qed

\end{document}